\documentclass[preprint]{elsarticle}
\usepackage[utf8]{inputenc}
\usepackage[T1]{fontenc}
\usepackage{graphicx}
\usepackage{grffile}
\usepackage{longtable}
\usepackage{wrapfig}
\usepackage{rotating}
\usepackage[normalem]{ulem}
\usepackage{amsmath}
\usepackage{textcomp}
\usepackage{amssymb}
\usepackage{capt-of}
\usepackage{hyperref}
\usepackage[disable]{todonotes}
\usepackage{cleveref}
\Crefname{equation}{}{}
\newcommand{\inlinelatex}[1]{#1}
\usepackage{subfigure}
\newcommand{\pvint}{\mathop{\mathrm PV}\!\int}
\newcommand{\pv}{\mathop{\mathrm PV}}
\usepackage{amsthm}
\theoremstyle{plain}
\newtheorem{theorem}{Theorem}
\newtheorem{heuristic}{Heuristic}

\newtheorem{corollary}[theorem]{Corollary}
\newtheorem{remark}[theorem]{Remark}
\newtheorem{definition}[theorem]{Definition}
\newtheorem{proposition}[theorem]{Proposition}
\date{\today}
\title{}
\hypersetup{
 pdfauthor={},
 pdftitle={},
 pdfkeywords={},
 pdfsubject={},
 pdfcreator={Emacs 26.3 (Org mode 9.3)}, 
 pdflang={English}}
\begin{document}

\listoftodos
\newpage

\begin{frontmatter}

\title{An Integral Equation Method for the Cahn-Hilliard Equation
       in the Wetting Problem}

\author[uiuc]{Xiaoyu Wei\corref{cor1}}
\ead{xywei@illinois.edu}

\author[njit]{Shidong Jiang\fnref{fn2}}
\ead{shidong.jiang@njit.edu}

\author[uiuc]{Andreas Kl{\"o}ckner\fnref{fn3}}
\ead{andreask@illinois.edu}

\author[hkust]{Xiao-Ping Wang\fnref{fn1}}
\ead{mawang@ust.hk}

\address[njit]{Department of Mathematical Sciences, New Jersey Institute
               of Technology, Newark, New Jersey 07102}
\address[uiuc]{Department of Computer Science, University of Illinois
               at Urbana-Champaign, Urbana, IL 61801}
\address[hkust]{Department of Mathematics, the Hong Kong University of
                Science and Technology, Clear Water Bay, Kowloon,
                Hong Kong, China}

\cortext[cor1]{Corresponding author}
\fntext[fn1]{This research was supported in part by the Hong Kong
             RGC-GRF grants 605513 and 16302715, RGC-CRF grant C6004-14G,
             and NSFC-RGC joint research grant N-HKUST620/15.}
\fntext[fn2]{This research was supported by the NSF under grant DMS-1418918.}
\fntext[fn3]{This research was supported by the NSF under grant DMS-1654756.}


\begin{abstract}
We present an integral equation approach to solving the Cahn-Hilliard equation
equipped with boundary conditions that model solid surfaces with prescribed
Young's angles. The discretization of the system in time using convex splitting leads to
a modified biharmonic equation at each time step.
To solve it, we split the solution into
a volume potential computed with free space kernels, plus
the solution to a second kind integral equation (SKIE).
The volume potential is evaluated with the help of a box-based
volume-FMM method. For non-box domains, source density is
extended by solving a biharmonic Dirichlet problem.
The near-singular boundary integrals are computed using
quadrature by expansion (QBX) with FMM acceleration.
Our method has linear complexity in the number of surface/volume degrees of freedom
and can achieve high order convergence with adaptive refinement
to manage error from function extension.
\end{abstract}

\end{frontmatter}

\section{Introduction}
\label{sec:orgf0fcb11}
\label{sec:introduction}

The Cahn-Hilliard equation is frequently used in the phase field model
in modeling processes involving the evolution of interfaces
(\cite{jacqmin_calculation_1999,liu_phase_2003,chen_phasefield_2002}).
There has been a large body of work on numerical methods for the Cahn-Hilliard
equations (cf.
\cite{barrett_finite_1999,demello_numerical_2005,feng_fourier_2009,shen_numerical_2010}
and references therein).
Most of the existing methods are finite element/finite difference methods,
or spectral methods for simple geometries.
Being a fourth order partial differential equation, standard finite difference or finite
element methods lead to very ill-conditioned matrices. This in turn causes low
accuracy and may require preconditioning for efficient solution.
Spectral methods are superior in terms of convergence order and number of
unknowns in the discretized system at the cost of only being able to
handle simple geometries.

In this paper, we present
a high-order accurate numerical method along with an efficient solution algorithm
for solving the two-phase
Cahn-Hilliard equation using integral equation methods in complex
geometry. At each time step,
our method requires the evaluation of
volume potentials on an adaptive volume grid, followed by the solution of
a well-conditioned system of second kind integral equations (SKIEs) with
unknown boundary densities. The volume potentials are
evaluated in linear complexity with respect to the number of quadrature nodes
using the volume fast multipole method (FMM)
based on a ``box FMM'' similar to, for example,
\cite{cheng_adaptive_2006,ethridge_new_2001,langston_freespace_2011,malhotra_pvfmm_2015,gholami_fft_2016}.
The boundary integral equations involve weakly singular
integrals, which are discretized using
the quadrature-by-expansion (QBX) method \cite{klockner_quadrature_2013}.
The idea of QBX is to exploit the
smoothness of the (smooth) layer potential away from the surface by forming
locally-valid expansions which are then evaluated to compute the
near or on-surface value of the potential.
Finally, the resulting linear system is solved by iterative
solvers such as GMRES \cite{saad_gmres_1986}
with matrix-vector products accelerated by a version of the FMM
\cite{greengard_fast_1987,cheng_fast_1999,gimbutas_generalized_2003,ying_kernelindependent_2004},
yielding linear complexity with respect to the number of boundary nodes.
When an adaptive volume mesh is used, the adaptive box structure is adjusted after
each step to capture the moving interface by refining in the neighborhood of the
interface and coarsening away from the interface based on solution gradient; in
addition, local refinement near the domain boundary is used to achieve higher
order accuracy during volume potential evaluation. The overall scheme is well
conditioned, potentially high-order in space, and has asymptotically optimal
complexity in terms of problem size.

The outline of the remainder of the paper is as follows.
In Section \ref{sec:problem-specification}, we present the model equations and
the temporal discretization scheme.
In Section \ref{sec:fundamental-solutions}, we analyze the static problem at each time step,
and provide a mathematical basis for the SKIE formulation.
In Section \ref{sec:integral-equation-formulation}, we present the SKIE
formulation for the static problem.
In Section \ref{sec:numerical-algorithms}, we present the numerical algorithms
to accompany the method formulation so that the final solver is
robust and has linear complexity, and in Section \ref{sec:numerical-results}
we present some numerical results.
In Section \ref{sec:conclusions}, we summarize the paper and discuss some future
perspectives.

\section{The Cahn-Hilliard Equation: Original Problem Specification}
\label{sec:orga9d9bd3}
\label{sec:problem-specification}

We consider a two-phase Cahn-Hilliard equation derived from
a simple phase-field model given by the following Ginzburg-Landau energy
functional
\cite{cahn_free_1958}:
\begin{equation}
\label{eqn:Ginzburg-Landau-energy}
\mathcal{E}_{GL}[\phi]
 = \int_\Omega \frac{ \epsilon }{ 2 } \| \nabla\phi(\mathbf{x}) \|^2
+ \frac{ [\phi(\mathbf{x})^2-1]^2 }{ 4 \epsilon } d\mathbf{x},
\end{equation}
where \(\phi(\mathbf{x})\) is the composition field
and \(\epsilon\) represents the interface thickness.
When \(\Omega=\mathbb{R}^2\), taking the gradient flow of
 (\ref{eqn:Ginzburg-Landau-energy})
in \(H^{-1}(\mathbb{R}^2)\)
yields the Cahn-Hilliard equation
\cite{cahn_spinodal_1971}:
\begin{align}
    \frac{\partial \phi(\mathbf{x},t)}{\partial t} & =
    \Delta \mu(\mathbf{x},t)
    \quad & \forall (\mathbf{x},t) \in \Omega\times\mathbb{R}^+
    \label{eqn:Cahn-Hilliard-phi},
    \\
    \mu(\mathbf{x},t) & = -\epsilon \Delta \phi(\mathbf{x},t)
      + \frac{ \phi(\mathbf{x},t)^3 - \phi(\mathbf{x},t) }{ \epsilon }
    \quad & \forall (\mathbf{x},t) \in \Omega\times\mathbb{R}^+
    \label{eqn:Cahn-Hilliard-mu},
\end{align}
where \(\mu\) is the \inlinelatex{Fr\'echet} derivative of
(\ref{eqn:Ginzburg-Landau-energy}),
also known as the chemical potential.

In practice, we are interested in solving the initial-boundary value problem
on bounded domains.
The initial conditions are
\begin{align}
  \phi(\mathbf{x},0)& = \phi_0(\mathbf{x})
  \quad & \forall \mathbf{x} \in \Omega
  , \label{eqn:cahn-hilliard-init-phi} \\
  \mu (\mathbf{x},0)& = -\epsilon \Delta \phi_0(\mathbf{x})
   + \frac{ \phi_0^3(\mathbf{x}) - \phi_0 (\mathbf{x})}{ \epsilon }
  \quad & \forall \mathbf{x} \in \Omega
  .
  \label{eqn:cahn-hilliard-init_mu}
\end{align}

Being a fourth order PDE, the Cahn-Hilliard equation requires two sets of
boundary conditions. In the wetting problem, the boundary models solid surfaces,
and a physically relevant way to impose boundary conditions is by adding a
surface energy term to the free energy to account for the interaction with the
solid wall:
\cite{modica_gradient_1987}
\begin{equation}
\label{eqn:solid-surface-energy}
    \mathcal{E}_S = \int_{\partial\Omega} \gamma(\phi(\mathbf{x})) ds_\mathbf{x},
\end{equation}
where we choose \(\gamma(\phi) = \frac{ \sqrt{2} }{ 3 } \cos \theta_Y \sin
\left( \frac{ \pi }{ 2 } \phi \right)\), which
gives the equilibrium contact angle (Young's angle)
of \(\theta_Y\)
\cite{xu_derivation_2010}.
Re-take the gradient flow but now with the total free energy
\(\mathcal{F} = \mathcal{E}_{GL} + \mathcal{E}_S\),
and the extra boundary energy term yields the relaxation boundary condition
\cite{chen_analysis_2014}
\begin{align}
    \frac{\partial \phi(\mathbf{x},t)}{\partial t} & =
    -\epsilon \partial _n\phi(\mathbf{x},t)
    + \frac{\partial \gamma(\phi)}{\partial \phi}(\mathbf{x},t)
    \quad & \forall (\mathbf{x},t) \in \partial\Omega\times\mathbb{R}^+
    \label{eqn:cahn-hilliard-solid-bc1},
\end{align}
as well as the zero-flux boundary condition
\begin{align}
    \partial _n \mu(\mathbf{x},t) & = 0
    \quad & \forall (\mathbf{x},t) \in \partial\Omega\times\mathbb{R}^+
    \label{cahn-hilliard-solid-bc2},
\end{align}
where \(\partial _n = \mathbf{n} \cdot \nabla\) with \(\mathbf{n}\) being the
 unit outward normal vector of \(\partial\Omega\).

Using Rothe's method, an energy decaying time stepping scheme for the system
(\ref{eqn:Cahn-Hilliard-phi} -- \ref{eqn:cahn-hilliard-init_mu}) and
(\ref{eqn:cahn-hilliard-solid-bc1}) -- (\ref{cahn-hilliard-solid-bc2})
is introduced in
\cite{gao_gradient_2012},
which is based on the idea of convex splitting
\cite{eyre_unconditionally_1997}:
at the \emph{n}-th time step, given \(\phi^n(\mathbf{x})\) ,
find \(\phi^{n+1}(\mathbf{x})\) such that
for all \(\mathbf{x}\) inside the domain,
\begin{align}
 \frac{ \phi^{n+1}(\mathbf{x}) - \phi^{n}(\mathbf{x})}{ \delta t }
  & = \Delta \mu ^{n+1} (\mathbf{x})
  , \label{eqn:convex-splitting-phi}\\
 \mu^{n+1} (\mathbf{x}) & = -\epsilon \Delta \phi^{n+1}(\mathbf{x}) +
  \frac{ s \phi^{n+1}(\mathbf{x}) - (1+s)\phi^n(\mathbf{x}) +
  (\phi^n)^3(\mathbf{x}) }{ \epsilon }
  , \label{eqn:convex-splitting-mu}
\end{align}
and that for all \(\mathbf{x}\) on the boundary,
\begin{align}
 \frac{ \phi^{n+1}(\mathbf{x}) - \phi^{n}(\mathbf{x})}{ \delta t }
  & = -\epsilon \partial_n\phi^{n+1}(\mathbf{x}) +
  \frac{\partial \gamma}{\partial\phi }(\phi^n(\mathbf{x}))
  , \label{eqn:convex-splitting-bc-phi} \\
 \partial_n \mu^{n+1} (\mathbf{x}) & = 0. \label{eqn:convex-splitting-bc-mu}
\end{align}
The splitting parameter \(s\) is chosen to be large enough such that the scheme is
energy-stable (under the assumption that \(\Vert\phi\Vert_\infty <
\infty\), see \cite{gao_gradient_2012} for details). If a bound to
the solution is known \(\Vert \phi \Vert_{\infty,
\Omega} \leq M\), \(\forall t > 0\), then letting \(s \geq (3M^2 - 1) / 2\)
guarantees energy stability of the scheme; however, even though the numerical
solutions mostly stay around the \([-1, 1]\) range,
the mathematical problem of solution boundedness for Cahn-Hilliard equation is
still open. In this paper we choose to fix \(s = 1.5\), which is sufficient to
allow for arbitrary \(\delta t\) in all our tests.

It is worth noting that time
discretization of the Cahn-Hilliard equation is still an active research
area. In this paper, we choose to adopt the first order linearly-implicit
convex-splitting method and direct our main efforts to
solving the static forced modified biharmonic subproblems using integral
equation approach. Other time discretization methods include: higher-order
convex splitting methods
(\cite{baskaran_convergence_2013,shen_secondorder_2012})
, stabilization methods
(\cite{shen_numerical_2010,zhu_coarsening_1999})
, the method of invariant energy quadratization (IEQ)
(\cite{guillen-gonzalez_linear_2013,yang_linear_2016}),
and the recent scalar auxiliary variable (SAV) method
\cite{shen_scalar_2018}.
Some of those advanced time integration schemes (like the SAV method)
produce fourth-order static boundary value problems with
constant coefficients, similar to the one treated in this paper.
The integral equation formulation and fast algorithms in this
paper can thus be straightforwardly generalized to adopt such
time discretization techniques.

For convenience of deriving integral representations of the solution, we
rewrite the problem
(\ref{eqn:convex-splitting-phi})-(\ref{eqn:convex-splitting-bc-mu})
by collecting unknowns onto one side of the equation:
\begin{align}
 (\Delta^2 - b \Delta + c)\phi^{n+1}(\mathbf{x})& = f_1(\mathbf{x})
  \quad &(\mathbf{x}\in\Omega), \label{eqn:cahn-hilliard-rewrite-1} \\
 (\Delta - b) \phi^{n+1}(\mathbf{x})
  + \frac{ 1 }{ \epsilon }\mu^{n+1}(\mathbf{x}) &= f_2(\mathbf{x})
  \quad &(\mathbf{x}\in\Omega), \label{eqn:cahn-hilliard-rewrite-2} \\
 (\partial_n + c) \phi^{n+1} (\mathbf{x}) & = h(\mathbf{x})
  \quad &(\mathbf{x}\in\partial\Omega), \label{eqn:cahn-hilliard-rewrite-3} \\
 \frac{ 1 }{ \epsilon }  \partial_n \mu ^{n+1}(\mathbf{x}) & = 0
  \quad &(\mathbf{x}\in\partial\Omega), \label{eqn:cahn-hilliard-rewrite-4}
\end{align}
where \(b = \frac{ s }{ \epsilon^2 }, c = \frac{ 1 }{ \epsilon \delta t }\),
and the inhomogeneous terms \(f_1, f_2\) and the boundary data \(h\) are given
by
\begin{align}
 f_1(\mathbf{x}) & = f_1[\phi^n(\mathbf{x})]
                 = c \phi^n(\mathbf{x})  + \Delta f_2(\mathbf{x})
  , \label{eqn:cahn-hilliard-rewrite-f1} \\
 f_2(\mathbf{x}) & = f_2[\phi^n(\mathbf{x})]
                 = \frac{ (\phi^n(\mathbf{x}))^3 - (1+s)\phi^n(\mathbf{x}) }
                          { \epsilon^2 }
  , \label{eqn:cahn-hilliard-rewrite-f2} \\
 h(\mathbf{x}) & = h[\phi^n(\mathbf{x})]
               = c \phi^n(\mathbf{x}) + \frac{ 1 }{ \epsilon }
                   \frac{\partial \gamma}{\partial \phi } (\phi^n(\mathbf{x}))
  . \label{eqn:cahn-hilliard-rewrite-h}
\end{align}

\section{Fundamental Solutions and Their Properties}
\label{sec:org0d51a18}
\label{sec:fundamental-solutions}

Denote the roots of the quadratic equation \(x^2-bx+c=0\) by
\(\lambda_1^2\) and \(\lambda_2^2\).
Then the fourth order operator \(\Delta^2 - b\Delta + c\)
in (\ref{eqn:cahn-hilliard-rewrite-1}) can be factored as
\((\Delta - \lambda_1^2)(\Delta - \lambda_2^2)\).
To employ an IE method for the solution of (\ref{eqn:cahn-hilliard-rewrite-1}),
we seek a fundamental solution \(G_0\) that should solve, in a weak sense,
\begin{equation}
\label{eqn:cahn-hilliard-greensfunc-eqn}
(\Delta - \lambda_1^2)(\Delta - \lambda_2^2)G_0(\mathbf{x}, \mathbf{y})
= \delta(\mathbf{x} - \mathbf{y}).
\end{equation}
Motivated by the factorization structure, we recall the Green's function of
the Yukawa operator that satisfies (weakly)
\begin{equation}
\label{eqn:yukawa-greensfunc-eqn}
(\Delta - \lambda_i^2) G_i(\mathbf{x}, \mathbf{y})
= \delta(\mathbf{x} -  \mathbf{y})
\end{equation}
and is given by the expression
\begin{equation}
\label{eqn:yukawa-greensfunc-expression}
G_i(\mathbf{x}, \mathbf{y}) = - \frac{ 1 }{ 2\pi }
K_0(\lambda_i r), \quad i=1,2,
\end{equation}
where \(r = \|x-y\|_2\) and \(K_0\) is the
modified Bessel function of the second kind of order zero (see, for example,
\cite{abramowitz_handbook_1965}).
For (\ref{eqn:cahn-hilliard-greensfunc-eqn}),
we formulate the fundamental solution as follows:
\begin{definition}
The fundamental solution $G_0$ is given by the formula
\begin{equation}
G_0(\mathbf{x}, \mathbf{y}) =
 - \frac{ 1 }{ 2\pi } \frac{ 1 }{ \lambda_1^2 - \lambda_2^2 }
 [ K_0(\lambda_1 r) - K_0(\lambda_2 r) ].
\label{eqn:cahn-hilliard-greensfunc-G0}
\end{equation}
In the special case of $\lambda_1=\lambda_2=\lambda$,
the funcdamental solution $G_0$ is defined by taking the limit
$\lambda_1 \rightarrow \lambda_2$ of (\ref{eqn:cahn-hilliard-greensfunc-G0}),
and is given by the formula
\begin{equation}
G_0(x,y)=- \frac{1}{4\pi\lambda} rK_1(\lambda r),
\label{eqn:cahn-hilliard-greensfunc-G0-limit}
\end{equation}
where $K_1$ is the modified Bessel function of the second kind of order one.
\end{definition}

For simplicity of discussion, we assume that \(\lambda_1^2
\neq \lambda_2^2\) for the rest for this paper. In practice,
if \(\lambda_1^2 = \lambda_2^2\), the equality can always be broken
by slightly changing the values of the time discretization
parameters \(s\) and \(\delta t\).

\begin{proposition}
$G_0$ as defined by (\ref{eqn:cahn-hilliard-greensfunc-G0}) satisfies
(\ref{eqn:cahn-hilliard-greensfunc-eqn}).
\begin{proof}
By definition (\ref{eqn:cahn-hilliard-greensfunc-G0}),
\begin{align*}
G_0(\mathbf{x}, \mathbf{y})
& =
 - \frac{ 1 }{ 2\pi } \frac{ 1 }{ \lambda_1^2 - \lambda_2^2 }
 [ K_0(\lambda_1 r) - K_0(\lambda_2 r) ] \\
& =
 \frac{ 1 }{ \lambda_1^2 - \lambda_2^2 }
 \left[
 G_1(\mathbf{x}, \mathbf{y}) - G_2(\mathbf{x}, \mathbf{y})
 \right],
\end{align*}
therefore,
\begin{align*}
 \Delta G_0(\mathbf{x}, \mathbf{y}) =
 \frac{ 1 }{ \lambda_1^2 - \lambda_2^2 }
 \left[
\Delta G_1(\mathbf{x}, \mathbf{y}) - \Delta G_2(\mathbf{x}, \mathbf{y})
 \right].
\end{align*}
Using (\ref{eqn:yukawa-greensfunc-eqn}), we have
\begin{align*}
 \Delta G_0(\mathbf{x}, \mathbf{y}) & =
 \frac{ 1 }{ \lambda_1^2 - \lambda_2^2 }
 \left\{
 \left[
\delta(\mathbf{x} - \mathbf{y}) + \lambda_1^2 G_1(\mathbf{x}, \mathbf{y})
 \right]
 -
 \left[
 \delta(\mathbf{x} - \mathbf{y}) + \lambda_2^2 G_2(\mathbf{x}, \mathbf{y})
 \right]
 \right\} \\
& =
 \frac{ \lambda_1^2 }{ \lambda_1^2 - \lambda_2^2 }
 G_1(\mathbf{x}, \mathbf{y}) -
 \frac{ \lambda_2^2 }{ \lambda_1^2 - \lambda_2^2 }
 G_2(\mathbf{x}, \mathbf{y}) \\
& =
 \lambda_1^2 G_0(\mathbf{x}, \mathbf{y}) +
 G_2(\mathbf{x}, \mathbf{y}) \\
& =
 \lambda_2^2 G_0(\mathbf{x}, \mathbf{y}) +
 G_1(\mathbf{x}, \mathbf{y}),
\end{align*}
that is,
\begin{equation}
 \label{eqn:relation-G0-G1-G2}
 (\Delta - \lambda_1^2) G_0(\mathbf{x}, \mathbf{y}) = G_2(\mathbf{x}, \mathbf{y}),
 \quad \text{and }
 (\Delta - \lambda_2^2) G_0(\mathbf{x}, \mathbf{y}) = G_1(\mathbf{x}, \mathbf{y}).
\end{equation}
Combining (\ref{eqn:relation-G0-G1-G2}) and (\ref{eqn:yukawa-greensfunc-eqn}),
we have
\begin{align*}
 (\Delta - \lambda_2^2)
 (\Delta - \lambda_1^2)
 G_0(\mathbf{x}, \mathbf{y}) =
 (\Delta - \lambda_1^2)
 (\Delta - \lambda_2^2)
 G_0(\mathbf{x}, \mathbf{y}) =
\delta(\mathbf{x} - \mathbf{y}).
\end{align*}
\end{proof}
\end{proposition}

Now we are equipped to define layer and volume potential operators.
\begin{definition}[Single layer potentials]
    Given a density function
    $\sigma(\mathbf{x}) \in C (\partial \Omega)$,
 the single layer potential operators are
    \begin{equation}
        \label{eqn:cahn-hilliard-single-layer-pot}
        S_i[\sigma](\mathbf{x}) = \int_{\partial\Omega}
        G_i(\mathbf{x}, \mathbf{y})
        \sigma(\mathbf{y}) d s_{\mathbf{y}}
        \qquad
         (\mathbf{x}\in\mathbb{R}^2,
         i = 0, 1, 2).
    \end{equation}
\label{def:ch-single-layer-potentials}
\end{definition}
\begin{definition}[Volume potentials]
    Given a density function
    $f(\mathbf{x}) \in C (\Omega)$
    , the volume potential operators are
    \begin{equation}
        \label{eqn:cahn-hilliard-volume-pot}
        V_i[f](\mathbf{x}) = \int_{\Omega}
        G_i(\mathbf{x}, \mathbf{y})
        f(\mathbf{y}) d \mathbf{y}
        \qquad
         (\mathbf{x}\in\mathbb{R}^2,
         i = 0, 1, 2).
    \end{equation}
\label{def:ch-volume-potentials}
\end{definition}
Note that here we only require the density functions to be continuous.
In practice, higher regularity is desired for
use of high order discretization, and our solver in this paper builds on
\(C^1\) density functions.

Regarding the jump conditions of layer potentials across the domain boundary,
we have the following result:
\begin{theorem}[Jump relations]
The normal derivative of $S_0[\sigma]$ has no jump across
the boundary $\partial\Omega$. However, when $\mathbf{x}$
approaches a point $\mathbf{p}\in\partial\Omega$
nontangentially, the normal derivative of the single layer
potential operators $S_i[\sigma]$ ($i=1,2$) on $\partial
\Omega$ satisfies the following jump relation:
\begin{equation}
    \label{eqn:Si-jump-relation}
    \lim_{\mathbf{x}\rightarrow\mathbf{p}^\pm}
    \frac{\partial S_i[\sigma](\mathbf{x})}
        {\partial \mathbf{n}_{\mathbf{p}}}
    = \left( \pm \frac{ 1 }{ 2 } I + \partial_n S_{i} \right)
    [\sigma](\mathbf{p})
    \qquad (i=1,2),
\end{equation}
where $\mathbf{n}_\mathbf{p}$ is the
 unit outward normal of $\partial\Omega$ at $\mathbf{p}$, and
\begin{equation}
    \label{eqn:principal-value-SLP}
    \partial_n S_i [\sigma](\mathbf{p}) =
    \pvint_{\partial\Omega}
    \frac{\partial G_i(\mathbf{p}, \mathbf{y})}{\partial \mathbf{n}_{\mathbf{p}}}
    \sigma(\mathbf{y}) ds_{\mathbf{y}}
    \qquad (i=1,2),
\end{equation}
and $I$ is the identity operator. Here
$\mathbf{x}\rightarrow\mathbf{p}^\pm$ means that $\mathbf{x}$ approaches
$\mathbf{p}$ from the exterior $(+)$ or the interior $(-)$ of the domain,
respectively, and $\pv \cdot$ denotes the Cauchy principal value.
\begin{proof}
We expand $K_0$ at the origin using asymptotic formulae
\cite[(9.6.10)--(9.6.13)]{abramowitz_handbook_1965}.
When $z \rightarrow 0$,
\begin{equation}
\label{eqn:cahn-hilliard-K0-expansion}
\begin{aligned}
K_0(z) & = - \ln \frac{ z }{ 2 } \sum_{k=0}^{\infty}
         \frac{ z^{2k} }{ 4^k (k!)^2 } + \sum_{k=0}^{\infty}
         \frac{ \psi(k+1)z^{2k} }{ 4^k (k!)^2 } \\
       & = - \left( \ln \frac{ z }{ 2 } + \gamma \right)
         \left( 1 + \frac{ z^2 }{ 4 } + \frac{ z^4 }{ 64 }
                + \mathcal{O}(z^6) \right)
         + \frac{ z^2 }{ 4 } + \frac{ 3z^4 }{ 128 } + \mathcal{O}(z^6),
\end{aligned}
\end{equation}
where $\psi(x)$ is the {\it psi} ({\it digamma}) function defined as
\begin{equation}
\label{eqn:psi-function}
\psi(n)=-\gamma + \sum_{k=1}^{n-1}
\frac{ 1 }{ k } \qquad (n =1, 2,3,\dots),
\end{equation}
and $\gamma=0.5772156649\ldots$ is Euler's constant
\cite{abramowitz_handbook_1965}.

From (\ref{eqn:cahn-hilliard-K0-expansion})
it is clear that the leading-order non-smooth behavior results from
the logarithmic term.
Then (\ref{eqn:Si-jump-relation})
follows directly from the well-known jump condition for the normal
derivative of the single
layer potential with logarithmic kernel
(\cite[Lemma 3.30]{folland_introduction_1976}).
\end{proof}
\label{thm:ch-lp-jump-relation}
\end{theorem}

Using the compactness of integral operators with weakly singular kernels
\cite[Theorem 1.11]{colton_integral_2013},
we also have the following corollary regarding the compactness
of those single layer potential operators and their derivatives
\begin{corollary}[Compactness]
For a Lipschitz domain $\Omega$, all three single layer
potential operators are compact in $C(\Omega)$. Also, the principal value parts of
$\partial_n S_i[\sigma]$ ($i=0,1,2$) are all compact.
\label{cor:ch-single-lp-compactness}
\end{corollary}

\section{Integral Equation Formulation}
\label{sec:orgb93e29a}
\label{sec:integral-equation-formulation}

In this section, we present a second-kind integral equation formulation
 for the system
(\ref{eqn:cahn-hilliard-rewrite-1})--(\ref{eqn:cahn-hilliard-rewrite-h}).
First we introduce a construction of the volume potential that
reduces the problem to a pure (volume-homogeneous) boundary value problem (PBVP).
Then we present an integral equation formulation for the PBVP that only
requires solving a Fredholm integral equation of the second kind.

\subsection{Reduction to Pure Boundary Value Problem}
\label{sec:org2bcf108}
\label{sec:reduction-to-pure}

We choose to represent the solutions \(\phi^n\), \(\mu^n\) in terms of sums of
volume and layer potentials. We first introduce the volume potential in the
terms of \(\tilde{u}\) and \(\tilde{v}\):
\begin{align}
 \phi^{n+1}(\mathbf{x}) & = \tilde{u}(\mathbf{x}) + u(\mathbf{x})
  \quad \text{with }
  \tilde{u}(\mathbf{x}) = V[\phi^n]
  \quad & (\mathbf{x} \in \Omega)
  , \label{eqn:cahn-hilliard-reduction-phi} \\
 \frac{ 1 }{ \epsilon } \mu^{n+1}(\mathbf{x}) & = \tilde{v}(\mathbf{x})
                                                  + v(\mathbf{x})
  \quad \text{with }
 \tilde{v}(\mathbf{x}) = f_2(\mathbf{x}) - (\Delta - b)\tilde{u}(\mathbf{x})
  \quad & (\mathbf{x} \in \Omega)
  , \label{eqn:cahn-hilliard-reduction-mu}
\end{align}
where we require that \(V[\phi^n]\) satisfies
\begin{equation}
  (\Delta^2 - b \Delta + c)V[\phi^n](\mathbf{x})  = f_1(\mathbf{x}),
  \quad (\mathbf{x} \in \Omega), \label{eqn:volume-potential-condition}
\end{equation}
and \(V[\cdot]\) is a volume potential operator to be defined.

For the sole purpose of removing inhomogeneities
in (\ref{eqn:cahn-hilliard-rewrite-1}) and (\ref{eqn:cahn-hilliard-rewrite-2}),
one could simply use the
volume potential
\begin{equation}
\label{eqn:volume-potential-V0}
  V_0[f_1] :=
           \int_\Omega G_0(\mathbf{x}, \mathbf{y}) f_1(\mathbf{y}) d\mathbf{y}.
\end{equation}
In numerical computation, however, evaluating \(f_1\)
(cf. (\ref{eqn:cahn-hilliard-rewrite-f1}))
directly is not desirable since it
involves second order differentiation and thus requires \(\phi^n \in C^2\),
which can be overly restrictive. Applying
Green's second identity to (\ref{eqn:volume-potential-V0}) yields
\begin{align*}
  V_0[f_1] & =  V_0[c\phi^n] - V_0[\Delta f_2] \\
           & =  V_0[c\phi^n] -
           \int_\Omega \Delta_{\mathbf{y}} G_0(\mathbf{x}, \mathbf{y}) f_2(\mathbf{y}) d\mathbf{y} \\
           & \quad +
           \int_{\partial \Omega}
               \partial_\mathbf{n_{\mathbf{y}}} G_0(\mathbf{x}, \mathbf{y}) f_2(\mathbf{y})
                d\mathbf{y}
           -
           \int_{\partial \Omega}
                G_0(\mathbf{x}, \mathbf{y}) \partial_\mathbf{n_{\mathbf{y}}} f_2(\mathbf{y})
                d\mathbf{y} \\
           & =  V_0[c\phi^n] -
           \int_\Omega \Delta_{\mathbf{x}} G_0(\mathbf{x}, \mathbf{y}) f_2(\mathbf{y}) d\mathbf{y} \\
           & \quad +
           \int_{\partial \Omega}
               \partial_\mathbf{n_{\mathbf{y}}} G_0(\mathbf{x}, \mathbf{y}) f_2(\mathbf{y})
                d\mathbf{y}
           -
           \int_{\partial \Omega}
                G_0(\mathbf{x}, \mathbf{y}) \partial_\mathbf{n_{\mathbf{y}}} f_2(\mathbf{y})
                d\mathbf{y}.
\end{align*}
Our construction of \(V[\phi^n]\) from
(\ref{eqn:cahn-hilliard-reduction-phi})
is obtained by
dropping the last two boundary integral terms and making use of
(\ref{eqn:relation-G0-G1-G2}), leading us to define
\begin{equation}
\label{eqn:volume-potential-formulation}
  V[\phi^n] = c V_0[\phi^n] + V_1[f_2[\phi^n]] + \lambda_2^2 V_0[f_2[\phi^n]].
\end{equation}
Inserting
(\ref{eqn:volume-potential-formulation})
into
(\ref{eqn:volume-potential-condition})
shows that the volume condition is satisfied despite the dropping of
the boundary terms.
By using our construction
(\ref{eqn:volume-potential-formulation})
numerical differentiation is completely avoided in the volume potential
evaluation stage.

Clearly, one does not need to solve for but only \emph{evaluate}
\(\tilde{u}\) and \(\tilde{v}\).
Specifically,
given \(\phi^n\), we first evaluate \(\tilde{u}\) by computing the
volume integral \(V[\phi^n]\) directly with the help of the box FMM. In the same
FMM pass, we also obtain \(\Delta \tilde{u}\) by taking derivatives
of the local expansions.
The loss of accuracy incurred due to numerical differentiation
can be recovered by increasing the FMM order.
 Then we can compute \(\tilde{v}\) directly
from its defining formula (\ref{eqn:cahn-hilliard-reduction-mu}).

Then \(u\) and \(v\) are the solutions of the following pure boundary value problem
\begin{align}
 (\Delta^2 - b \Delta + c)u(\mathbf{x})  & = 0
  \quad & (\mathbf{x} \in \Omega), \label{eqn:cahn-hilliard-pure-u} \\
 v + (\Delta - b)u(\mathbf{x})  & = 0
  \quad & (\mathbf{x} \in \Omega), \label{eqn:cahn-hilliard-pure-v} \\
 (\partial_n + c) u(\mathbf{x}) &= g_1(\mathbf{x})
  \quad & (\mathbf{x} \in \partial\Omega),
  \label{eqn:cahn-hilliard-pure-bc-u} \\
 \partial_n  v(\mathbf{x}) &= g_2(\mathbf{x})
  \quad & (\mathbf{x} \in \partial\Omega),
  \label{eqn:cahn-hilliard-pure-bc-v}
\end{align}
where the boundary data \(g_1\) and \(g_2\) are given by
\begin{align}
g_1(\mathbf{x})  & = h(\mathbf{x}) - \tilde{u}_n(\mathbf{x})
                                   - c \tilde{u}(\mathbf{x}),
 \label{eqn:cahn-hilliard-pure-g1} \\
g_2(\mathbf{x})  & = - \tilde{v}_n(\mathbf{x}).
 \label{eqn:cahn-hilliard-pure-g2}
\end{align}

\subsection{Second Kind Integral Equation Formulation for the Pure Boundary Value Problem}
\label{sec:orge9fd3c0}
\label{sec:skie-for-pure}

We now derive a second kind integral equation (SKIE) formulation for the
boundary value problem
\Crefrange{eqn:cahn-hilliard-pure-u}{eqn:cahn-hilliard-pure-g2}.
We first represent \(u\) by the formula
\begin{equation}
\label{eqn:cahn-hilliard-representation-u}
    \begin{aligned}
        u(\mathbf{x}) & = S_1[\sigma_1](\mathbf{x})
                        + S_0[\sigma_2](\mathbf{x}) \\
                      & = \int_{\partial\Omega} \left[
                          G_1(\mathbf{x}, \mathbf{y}) \sigma_1(\mathbf{y})
                        + G_0(\mathbf{x}, \mathbf{y}) \sigma_2(\mathbf{y})
                          \right] d s_{\mathbf{y}},
    \end{aligned}
\end{equation}
where \(\sigma_i\) (\(i=1,2\)) are unknown densities on \(\partial\Omega\).
Obviously, this representation satisfies
(\ref{eqn:cahn-hilliard-pure-u}). Substituting
(\ref{eqn:cahn-hilliard-representation-u}) into
(\ref{eqn:cahn-hilliard-pure-v}) yields
\begin{equation}
\label{eqn:cahn-hilliard-representation-v}
    \begin{aligned}
        v(\mathbf{x}) & = (-\Delta + b) u(\mathbf{x}) \\
                      & = \int_{\partial\Omega} \left[
                        (-\Delta + b) G_1(\mathbf{x}, \mathbf{y}) \sigma_1(\mathbf{y})
                      + (-\Delta + b) G_0(\mathbf{x}, \mathbf{y}) \sigma_2(\mathbf{y})
                          \right] d s_{\mathbf{y}},
    \end{aligned}
\end{equation}

Note that since \(b = \lambda_1^2 + \lambda_2^2\), the following corollary
follows naturally from the definition of the Yukawa kernels
(\ref{eqn:yukawa-greensfunc-eqn}):
\begin{corollary}
    When $\mathbf{y} \neq \mathbf{x}$,
    $G_i(\mathbf{x}, \mathbf{y})$ ($i=1,2$) are both locally eigenfunctions of
    the Yukawa operator $\Delta-b$. Specifically,
    \begin{equation}
        \label{eqn:Gi-are-eigenfunctions}
        (\Delta-b) G_1(\mathbf{x}, \mathbf{y}) = -\lambda_2^2 G_1(\mathbf{x}, \mathbf{y}), \quad
        (\Delta-b) G_2(\mathbf{x}, \mathbf{y}) = -\lambda_1^2 G_2(\mathbf{x}, \mathbf{y}).
    \end{equation}
    \begin{proof}
        We only prove the first identity. The other one follows analogously.
        Substituting $b = \lambda_1^2 + \lambda_2^2$ into the left hand side
        of the identity,
        \begin{equation*}
        (\Delta-b) G_1(\mathbf{x}, \mathbf{y}) =
        [\Delta-(\lambda_1^2 + \lambda_2^2)] G_1(\mathbf{x}, \mathbf{y}) =
        (\Delta-\lambda_1^2) G_1(\mathbf{x}, \mathbf{y})
        - \lambda_2^2 G_1(\mathbf{x}, \mathbf{y}).
        \end{equation*}
        By definition, $(\Delta-\lambda_1^2) G_1(\mathbf{x}, \mathbf{y}) =
        \delta(\mathbf{x}, \mathbf{y}) = 0$. Thus
        \begin{equation*}
        (\Delta-b) G_1(\mathbf{x}, \mathbf{y}) = - \lambda_2^2 G_1(\mathbf{x}, \mathbf{y}).
        \end{equation*}
    \end{proof}
\label{cor:Gi-are-eigenfunctions}
\end{corollary}

On the other hand, (\ref{eqn:relation-G0-G1-G2}) yields
\begin{equation}
\label{eqn:simplify-laplace-G0}
    \begin{aligned}
        (-\Delta+b) G_0(\mathbf{x}, \mathbf{y})
        & = [-\Delta + (\lambda_1^2 + \lambda_2^2)] G_0(\mathbf{x}, \mathbf{y}) \\
        & = - G_1(\mathbf{x}, \mathbf{y}) + \lambda_1^2 G_0(\mathbf{x}, \mathbf{y}) \\
        & = - G_2(\mathbf{x}, \mathbf{y}) + \lambda_2^2 G_0(\mathbf{x}, \mathbf{y}).
    \end{aligned}
\end{equation}

Thus
\begin{equation}
\label{eqn:simplified-representation-v}
v(\mathbf{x}) = \int_{\partial\Omega}\left\{
\lambda_2^2 G_1(\mathbf{x}, \mathbf{y}) \sigma_1(\mathbf{y})
+ [-G_1(\mathbf{x}, \mathbf{y}) + \lambda_1^2 G_0(\mathbf{x}, \mathbf{y})]
 \sigma_2(\mathbf{y})
\right\} d s_{\mathbf{y}}.
\end{equation}

Combining the jump relations (\ref{eqn:Si-jump-relation}),
the boundary conditions
(\ref{eqn:cahn-hilliard-pure-bc-u}) and (\ref{eqn:cahn-hilliard-pure-bc-v}),
and the representations
(\ref{eqn:cahn-hilliard-representation-u}) and
(\ref{eqn:simplified-representation-v}), we obtain the following system of
boundary integral equations:

\begin{equation}
\label{eqn:cahn-hilliard-boundary-integral-equations}
    (D + A)[\sigma](\mathbf{x}) = g(\mathbf{x})
    \qquad (\mathbf{x}\in\partial\Omega),
\end{equation}
where
\begin{equation}
\label{eqn:cahn-hilliard-boundary-integral-equation-components}
    D = - \frac{ 1 }{ 2 } \begin{bmatrix}
    1 & 0 \\ \lambda_2^2 & - 1
    \end{bmatrix},
    \quad
    A = \begin{bmatrix}
    A_{11} & A_{12} \\ A_{21} & A_{22}
    \end{bmatrix},
    \quad
    \sigma = \begin{bmatrix}
    \sigma_{1} \\ \sigma_{2}
    \end{bmatrix},
    \quad
    g = \begin{bmatrix}
    g_{1} \\ g_{2}
    \end{bmatrix},
\end{equation}
and the entries of the operator matrix \(A\) are given by the formulae
\begin{equation}
\label{eqn:cahn-hilliard-boundary-integral-equation-A-entries}
\begin{aligned}
A_{11} &= \partial_n S_1 + c S_1, \quad
&A_{12} = \partial_n S_0 + c S_0, \\
A_{21} &= \lambda_2^2 \partial_n S_1, \quad
& A_{22} = - \partial_n S_1 + \lambda_1^2 \partial_n S_0.
\end{aligned}
\end{equation}
Clearly, \(D\) has nonzero determinant, and all entries of \(A\) are compact
operators (Corollary \ref{cor:ch-single-lp-compactness}). Therefore, the
system (\ref{eqn:cahn-hilliard-boundary-integral-equations}) is a Fredholm
integral equation of the second kind.
When discretizing a
second kind integral equation, the condition number of the
resulting linear system remains bounded when the number of
unknowns increases as the discretization gets
refined. For an iterative solver achieving a fixed residual norm,
the numerical solution thus remains accurate as the mesh is refined.
Besides, with bounded
condition number, the number of GMRES iterations is less likely
to grow when refining the mesh.

\section{Numerical Realization: Discretization and Algorithms}
\label{sec:org7e17165}
\label{sec:numerical-algorithms}

For the implementation of our method, the volume potentials in
(\ref{eqn:volume-potential-formulation})
are evaluated using a version of the the volume FMM
(\cite{cheng_adaptive_2006,ethridge_new_2001,langston_freespace_2011,malhotra_pvfmm_2015,gholami_fft_2016}).
We further discretize the SKIE system
(\ref{eqn:cahn-hilliard-boundary-integral-equations})
using the Nyström method using QBX quadrature and solve the linear system
with the help of GMRES. The required matrix-vector products are carried out
using the GIGAQBX quadrature-enabled fast algorithm
\cite{klockner_quadrature_2013,rachh_fast_2017,wala_fast_2018}.
To handle complex input geometries, we use
Gmsh \cite{geuzaine_gmsh_2009} to
generate a triangulation along the domain boundary with each
panel roughly of the size \(h_b\). Over each boundary panel,
Legendre-Gauss quadrature nodes of order \(q_b\) are used to form
the discretization.

If \(\Omega = [x_l, x_r]^2\), the volume FMM bridges seamlessly with the boundary
discretization; however, additional treatment is needed when dealing with
general geometries, as the volume FMM requires box-shaped partitions to operate.
In this section, we present those extra steps necessary to achieve flexibility
in handling complex geometries.

\subsection{Regularity and Non-Square Domains}
\label{sec:org3122c74}
For non-square domains, the volume potentials are evaluated over a box
\(B\) that encloses the physical domain such that \(\Omega \subset B\).
Use of the box FMM requires that
volume source densities be extended to \(B\setminus\Omega\).
Noting that \(f_1, f_2\) are all algebraic functions of \(\phi^n\), we only need
to extend \(\phi\) to \(B\setminus \Omega\).
There are several
factors that need to be considered when deciding how to extend \(\phi\):
\begin{itemize}
\item The smoother the extension is, the higher accuracy in the order of
      approximation for
      those densities can be achieved using (piecewise) polynomials, yielding
      better overall accuracy.
\item Higher order extensions are more costly to evaluate. Typically for
      \(C^k\) extension one needs to solve a PDE where \(k+1\) boundary conditions
      can be imposed, e.g., \(\Delta^{k+1} u = 0\). It is natural to ask that the
      extension process does not cost much more than solving the original
      problem, i.e., \(k\leq1\).
\item Our SKIE formulation poses certain minimum regularity requirements.
      Specifically, evaluation of \(g_2\) in (\ref{eqn:cahn-hilliard-pure-bc-v})
      requires at minimum \(\tilde{v}\in C^1(B)\), which then requires \(\tilde{u}\)
      to be in \(C^3(B)\).
\end{itemize}

To satisfy the condition \(\tilde{u}\in C^3(B)\), the density \(\phi\) should
satisfy certain regularity requirements.
First we note the standard regularity result for the Poisson equation
\cite{gilbarg_elliptic_2001}:
\begin{theorem}
\label{thm:poisson-vp-estimate}
Let $\Omega \subset \mathbb{R}^d$ be open and bounded,
    \begin{equation}
    u(x) := \int_\Omega \Phi(\mathbf{x} - \mathbf{y}) f(\mathbf{y})
    d \mathbf{y},
    \end{equation}
where $\Phi$ is the fundamental solution. Then if
$f \in C_0^\alpha(\overline{\Omega})$, $0 < \alpha < 1$, then
$u \in C^{2,\alpha}(\overline{\Omega})$, and
    \begin{equation}
    \lVert u(x) \rVert_{C^{2,\alpha}(\overline{\Omega})} \leq
    c \lVert f \rVert _{C^\alpha(\overline{\Omega})}.
    \end{equation}
\end{theorem}

The leading order singular (least regular) term in our \(V[\phi^n]\)
construction is \(V_1[f_2]\), for which we have the following result:
\begin{proposition}
 If $f_2 \in C^1(B)$, the volume potential $V_1[f_2] \in C^3(B)$.
\begin{proof}
Since, like the Poisson kernel, the kernel $G_1$ also has a leading-order
logarithmic singularity,
Theorem \ref{thm:poisson-vp-estimate} also applies to $V_1[f_2]$.
Combined with the assumption that the source density is in $C^1(B)$,
the proposition follows.
\end{proof}
\end{proposition}

Based on these considerations, we choose to use \(C^1\) extensions.
To obtain such an extension, we assume that \(\Omega\) is simply
connected and solve the following Stokes problem:
\begin{align}
    -\Delta \mathbf{u} + \nabla p &= 0
    & \text{in } \mathbb{R}^2\setminus \Omega, \label{eqn:stokes-momentum}\\
    \nabla \cdot \mathbf{u} &= 0
    & \text{in } \mathbb{R}^2\setminus \Omega, \label{eqn:stokes-divergence-free}\\
    \mathbf{u} &= \nabla^\perp \phi
    & \text{on } \partial\Omega\label{eqn:stokes-bc},
\end{align}
where \(\mathbf{u}\) is the velocity field of the fluid and \(p\) is the pressure
field. The Stokes problem
(\ref{eqn:stokes-momentum}, \ref{eqn:stokes-divergence-free}, \ref{eqn:stokes-bc})
can also be solved via an SKIE formulation
\cite{greengard_integral_1996}.
We evaluate the stream function \(\omega\) up to a constant by evaluating complex
layer potential representations in \cite{rachh_integral_2018}, such that
\begin{equation}
\label{eqn:stream-function}
  \mathbf{u} = \nabla^\perp \omega =
  \begin{bmatrix}
  \frac{\partial \omega}{\partial x_2} \\
  - \frac{\partial \omega}{\partial x_1} \\
  \end{bmatrix},
\end{equation}
then \(\omega\) is a biharmonic function that solves the boundary value problem
\begin{align}
    \Delta^2\omega &= 0,
    & \text{in } \mathbb{R}^2\setminus \Omega, \label{eqn:biharmonic-streamfunc}\\
    \nabla \omega &= \nabla \phi,
    & \text{on } \partial\Omega\label{eqn:biharmonic-streamfunc-bc}.
\end{align}
Since \(\partial\Omega\) is simply connected, we only need to add a constant to
\(\omega\) to obtain a \(C^1\) extension of \(\phi\):
\begin{equation}
\label{eqn:phi-extension}
  \phi(\mathbf{x}) = \omega(\mathbf{x}) + \left(
    \int_{\partial\Omega} \phi(\mathbf{y})  ds_\mathbf{y} -
    \int_{\partial\Omega} \omega(\mathbf{y}) ds_\mathbf{y}
  \right)
\qquad
(\mathbf{x} \in \mathbb{R}^2\setminus \Omega)
.
\end{equation}

\begin{remark}
The assumption of \(\Omega\) being simply connected is for simplicity of
the formulation and not critical.
If \(\Omega\) is topologically more complicated, methods in
 \cite{farkas_mathematical_1989,rachh_integral_2018}
can still be used for performing \(C^1\) extension.
\end{remark}

\subsection{Boundary Layers}
\label{sec:org3ac9ad3}

From matched asymptotic analysis of the Cahn-Hilliard equation
(\ref{eqn:Cahn-Hilliard-phi})--(\ref{cahn-hilliard-solid-bc2})
 \cite{chen_analysis_2014a}, it is apparent that for very small \(\epsilon\),
the solution can develop boundary layers.

Boundary layers cause stability concerns for our scheme using
\(C^1\) extension. Since the exterior biharmonic problem does not
admit a Laplace-type maximum principle, the presence of the
boundary layer causes the extended density to be very large in
magnitude compared to the density inside the domain, resulting in
lost accuracy after each time step. As a result, the numerical
solution may not remain bounded after several time steps and lose
energy stability.

Dealing with such boundary layers numerically is the subject of ongoing research.
Using an adaptive volume mesh that refines towards the boundary can be helpful.
Instead of evaluating the volume potentials over one single bounding box \(B\), we
can choose \(B\) to be union of a set of boxes that covers \(\Omega\). By reducing
the Hausdorff distance between \(\partial B\) and \(\partial\Omega\), the
extended density is less likely to blow up.
However, such refinement is very expensive in practice, and there
is no guarantee that the boundary layer is resolved a-priori
since a reliable refinement criterion is not known.

For long-time simulations, we provide an alternative formulation which is
a less costly approximation based on the total energy of the system.
When solving the SKIE
(\ref{eqn:cahn-hilliard-boundary-integral-equations}), we
use a homogeneous right hand side for the second equation \(g_2 = 0\),
which effectively makes the layer potential \(u\) satisfy the contact angle
dynamics, while ignoring the second boundary condition.
After solving for \(\sigma_1\) and evaluating \(u\), we add a constant
to \(u\) such that the integral of \(\phi\) over the whole domain is
conserved.
This approach conserves the total mass and keeps the numerical
solution bounded at all time.

Formally, the approximation replaces the representation
(\ref{eqn:cahn-hilliard-representation-u})
with the following stabilized representation:
\begin{equation}
\label{eqn:cahn-hilliard-representation-u-stabilized}
    \begin{aligned}
        u(\mathbf{x}) & = S_1[\sigma_1](\mathbf{x}) + C
    \end{aligned}
\end{equation}
where \(C\) is the added constant. The net effect of the approximation spreads the
influence of boundary layers
across the whole domain. Long-time simulations using the stabilized
representation are performed and the results are in line with those obtained
from finite element method. Although the modification seems crude, it maintains
the most relevant solution characteristics including:

\begin{itemize}
  \item $H^{-1}$ gradient flow structure in the bulk region, since the volume potential is unchanged.
  \item Contact angle dynamics from the first boundary condition.
  \item Integral conservation over the whole domain.
\end{itemize}

\section{Numerical Results}
\label{sec:org6133647}
\label{sec:numerical-results}
\subsection{Spatial Convergence}
\label{sec:orgff4d44d}
\label{sec:convergence-test}
Based on error analysis of the volume FMM (\cite{ethridge_new_2001}),
QBX (\cite{epstein_convergence_2013}),
 the GIGAQBX fast algorithm (\cite{klockner_quadrature_2013}, \cite{wala_fast_2018}) and
standard approximation theory (\cite{trefethen_approximation_2013}),
we expect the following error estimate for the overall numerical accuracy of
our scheme:
\begin{heuristic}[Numerical Accuracy]
Assuming density extension has regularity $\phi^0 \in \mathcal{C}^k(B)$,
 denote the numerical solution at $t = \delta t$ by $\phi_h^1$ and
 its relative error by
$E_0 := \lVert \phi_h^1 - \phi(\delta t) \rVert_\infty /
           \lVert \phi(\delta t) \rVert_\infty$.
Asymptotically, for a single time step, we have
\begin{equation*}
E_0 \leq C_0 \delta t + \left( E_v + E_b \right),
\end{equation*}
where $C_0$ is a constant independent of $\delta_t$,
$E_v$ includes the errors from volume potential evaluation
and $E_b$ includes the errors from the BIE solution.
\begin{enumerate}
\item
For $E_v$, we have
\begin{align*}
E_v \leq \underbrace{C_1 h_v^{\min(\sqrt{q_v}, k+1)}}_{
                  \text{discretization error}}
             +
      \underbrace{C_2 2^{-m_v}}_{
                  \text{volume FMM error}}
             + \xi,
\end{align*}
where $h_v$ is the volume mesh size,
$q_v$ is the number of quadrature points per volume cell,
$m_v$ is the volume FMM order,
$\xi$ is the error of near-field volume potentials (controlled by the tolerance of adaptive quadrature over nearfield boxes),
and $C_1, C_2$ are constants independent of $h_v$ and $\xi$.

\item
For $E_b$,
assuming that the solution density $\sigma$ is smooth and that
its Hölder norms $\lVert \sigma \rVert_{C^{s, \beta}}$ for
$s \in \mathbb{N}, \beta \in [0, 1)$
 can be controlled with $\lVert \phi(\delta t)\rVert_\infty$,
we have
\begin{align*}
E_b \leq & \underbrace{C_3 h_b ^ {q_b}}_{
                  \text{discretization error}}
      +
      \underbrace{C_4 2^{-2q_b}}_{
                  \text{QBX quadrature error}}
      +
      \underbrace{C_5 h_b ^ {p + 1} 2^{-p}}_{
                  \text{QBX truncation error}} \\
&      +
      \underbrace{C_6 2^{-m_b}}_{
                  \text{QBX FMM error}}
      + C_7 \eta,
\end{align*}
where $h_b$ is the boundary mesh size,
$q_b$ is the number of quadrature points per boundary cell,
$p$ is the QBX order,
$m_b$ is the FMM order used in GIGAQBX,
$\eta$ is the error of the linear solver (controlled by the tolerance of GMRES),
and $C_3, C_4, C_5, C_6, C_7$ are constants independent of $h_b$ and $\eta$.
\end{enumerate}
\label{thm:error-analysis}
\end{heuristic}

Among all the error sources listed in Heuristic \ref{thm:error-analysis},
we are most interested in the discretization errors that decay algebraically
with mesh refinement.
In this section, we describe experiments to verify this expected error
decay, by choosing other parameters of the scheme so as to minimize the contributions
of their corresponding errors.
Specifically, we fix the tolerance for GMRES to be \(10^{-14}\), and the tolerance for
(precomputed) near-field adaptive quadrature in the volume FMM to be \(10^{-13}\). Also, we set
the FMM orders and the QBX order to be sufficiently high.
Specifically, we let the QBX order match the boundary
quadrature order \(q_b = p + 1\),
and fix the volume FMM order and the GIGAQBX's FMM order \(m_v = m_b = 8\).
We use two-dimensional Taylor expansions for the FMM, thus achieving FMM
order \(8\) in \(2\)D requires \(\frac{1}{2} \times 8 \times 9 = 36\) terms.
To back our choices, neither increasing the QBX order \(p\) up to \(10\) nor
the FMM order \(m_v, m_b\) up to \(10\) yields obvious effects on the results,
indicating that the error in the computation of the volume potential
is dominated by discretization errors.
The same characteristic mesh size \(h_b = h_v = \delta x\) is used for volume and
boundary meshes. For simplicity, we also keep the orders of volume and boundary discretizations
the same.

To gauge convergence, we fix \(\delta t = 10^{-2}\) and perform self-convergence tests by refining
both spatial grids successively and calculating
the relative error from the result of one time step
\(e_0 := \lVert \phi^1_h - \phi^1_0 \rVert_\infty / \lVert\phi^1_h\rVert_\infty\),
where \(\phi^1_0\) is the reference solution.
The tests are performed on a circular domain with radius
\(0.247\) centered at the origin. The bounding box is chosen to
be \(B = [-0.25, 0.25]^2\). The initial condition is
\begin{equation}
\label{convergence-test-init}
  \phi^0(\mathbf{x}) = \tanh \left[
           \frac{10}{\sqrt{2\epsilon}}
           (\left| \mathbf{x} \cdot \mathbf{e}_1 \right| - 0.1)
           \right],
\end{equation}
where \(\epsilon=10^{-2}\), and \(\mathbf{e}_1 = (1, 0)\).

In this work, since we are using \(C^1\) extension of \(\phi\), based on the
estimates above, when the discretization error dominates,
the total error of the scheme should scale as
\(\mathcal{O}(h_b^{q_b}) + \mathcal{O}(h_v^{\min(\sqrt{q_v},
2)})\) when \(h_b, h_v \rightarrow 0\).
In particular, when \(h_b = h_v = \delta x\), and \(q_b = \sqrt{q_v} = r\),
the total discretization error is expected to be
\(\mathcal{O}(\delta x ^r)\) as \(\delta x \rightarrow 0\).

When using piecewise constant approximations for both boundary and
volume discretizations, the results are shown in Table
\ref{tab:convergence-1st-order},
where \(q_b, q_v\) stand for the number of quadrature points per cell on
the boundary mesh and the volume mesh, respectively.
The convergence order is calculated using values
of the error measure \(e_0\)
from the current row and the row above.
The same statistics for the second order discretizations (i.e.
piecewise linear approximations on the boundary mesh and piecewise
bi-linear approximations on the volume mesh)
are shown in Table
\ref{tab:convergence-2nd-order}.

In both cases, the numerical results align well with expectations.
When \(r=1\), \(q_b = q_v = 1\), discretizations for both boundary and volume
densities are piecewise constant (first order), and Table
\ref{tab:convergence-1st-order}
demonstrates empirical convergence order close to \(1\).
When \(r=2\), \(q_b = 2, q_v = 4\), boundary densities are discretized with
piecewise linear functions, and volume densities are discretized with
piecewise bilinear functions, in which case Table \ref{tab:convergence-2nd-order}
shows empirical convergence order close to \(2\).

When higher-order volume approximations are used with a uniform
volume mesh, the convergence order will be limited to
 \(2\) due to the volume densities being limited to only \(C^1\) at
the boundary. For example, as shown in Table \ref{tab:convergence-3rd-order-raw}.

\begin{table}[htbp]
\centering
\begin{tabular}{|c|c|c|c|c|}
\hline
\(q_b\) & \(q_v\) & \(\delta x\) & \(e_0\) & EOC\\
\hline
1 & 1 & \(1.56 \times 10^{-2}\) & \(9.91 \times 10^{-1}\) & --\\
1 & 1 & \(7.81 \times 10^{-3}\) & \(6.72 \times 10^{-1}\) & 0.56\\
1 & 1 & \(3.91 \times 10^{-3}\) & \(2.91 \times 10^{-1}\) & 1.21\\
1 & 1 & \(1.95 \times 10^{-3}\) & \(1.34 \times 10^{-1}\) & 1.19\\
1 & 1 & \(9.77 \times 10^{-4}\) & \(9.02 \times 10^{-2}\) & 0.57\\
1 & 1 & \(4.88 \times 10^{-4}\) & \(5.24 \times 10^{-2}\) & 0.78\\
\hline
\end{tabular}
\caption{\label{tab:convergence-1st-order}Results Using First-Order Approximation (\(q_b\): number of quadrature points per boundary cell, \(q_v\): number of quadrature points per volume cell, \(e_0\): relative error at \(t=\delta t\)).}

\end{table}

\begin{table}[htbp]
\centering
\begin{tabular}{|c|c|c|c|c|}
\hline
\(q_b\) & \(q_v\) & \(\delta x\) & \(e_0\) & EOC\\
\hline
2 & 4 & \(1.56 \times 10^{-2}\) & \(4.01 \times 10^{-1}\) & --\\
2 & 4 & \(7.81 \times 10^{-3}\) & \(2.41 \times 10^{-1}\) & 0.73\\
2 & 4 & \(3.91 \times 10^{-3}\) & \(6.99 \times 10^{-2}\) & 1.79\\
2 & 4 & \(1.95 \times 10^{-3}\) & \(1.76 \times 10^{-2}\) & 1.99\\
2 & 4 & \(9.77 \times 10^{-4}\) & \(4.84 \times 10^{-3}\) & 1.86\\
\hline
\end{tabular}
\caption{\label{tab:convergence-2nd-order}Results Using Second-Order Approximation (\(q_b\): number of quadrature points per boundary cell, \(q_v\): number of quadrature points per volume cell, \(e_0\): relative error at \(t=\delta t\)).}

\end{table}

However, higher order approximations can still be used by
adding extra levels of refinement
to volume cells that intersect \(\partial \Omega\).
To guide the mesh refinement process, we form a
two-dimensional tensor-product
Legendre series of degree \(q_v\) (using the discrete tensor-product
Legendre transform) over each cell, and compute the sum of the
absolute values of the \(q_v + 1\) coefficients of
the leading terms (terms of total degree \(q_v\)) as the error
indicator for that cell.
To make the adaptive mesh for the box FMM, we start from a
uniform mesh with prescribed
\(\delta x\), and iteratively perform refinement cycles as
follows:
\begin{enumerate}
\item Compute the average values of the
error indicators for the interior cells as $M_i$, and that for the
boundary-intersecting cells as $M_b$.
\item Refine each boundary-intersecting cell by one level (into four cells).
\item Refine some additional cells as necessary to keep the tree balanced.
\end{enumerate}
We repeat the refinement cycles until \(M_b \leq M_i\).
As shown by the results in Table \ref{tab:convergence-3rd-order},
since the source density is smooth away from the boundary,
by performing such refinements,
the solver is able to take advantage of higher order of convergence
in the bulk region where the source density is smooth.

\begin{table}[htbp]
\centering
\begin{tabular}{|c|c|c|c|c|}
\hline
\(q_b\) & \(q_v\) & \(\delta x\) & \(e_0\) & EOC\\
\hline
3 & 9 & \(1.56 \times 10^{-2}\) & \(2.11 \times 10^{-1}\) & --\\
3 & 9 & \(7.81 \times 10^{-3}\) & \(3.73 \times 10^{-2}\) & 2.50\\
3 & 9 & \(3.91 \times 10^{-3}\) & \(1.04 \times 10^{-2}\) & 1.84\\
3 & 9 & \(1.95 \times 10^{-3}\) & \(2.36 \times 10^{-3}\) & 2.13\\
3 & 9 & \(9.77 \times 10^{-4}\) & \(5.58 \times 10^{-4}\) & 2.08\\
\hline
\end{tabular}
\caption{\label{tab:convergence-3rd-order-raw}Results Using Third-Order Approximation, without Extra Boundary Refinement (\(q_b\): number of quadrature points per boundary cell, \(q_v\): number of quadrature points per volume cell, \(e_0\): relative error at \(t=\delta t\)).}

\end{table}

\begin{table}[htbp]
\centering
\begin{tabular}{|c|c|c|c|c|}
\hline
\(q_b\) & \(q_v\) & \(\delta x\) & \(e_0\) & EOC\\
\hline
3 & 9 & \(1.56 \times 10^{-2}\) & \(2.11 \times 10^{-1}\) & --\\
3 & 9 & \(7.81 \times 10^{-3}\) & \(3.45 \times 10^{-2}\) & 2.61\\
3 & 9 & \(3.91 \times 10^{-3}\) & \(4.42 \times 10^{-3}\) & 2.96\\
3 & 9 & \(1.95 \times 10^{-3}\) & \(6.99 \times 10^{-4}\) & 2.66\\
3 & 9 & \(9.77 \times 10^{-4}\) & \(1.04 \times 10^{-4}\) & 2.75\\
\hline
\end{tabular}
\caption{\label{tab:convergence-3rd-order}Results Using Third-Order Approximation, with Extra Boundary Refinement (\(q_b\): number of quadrature points per boundary cell, \(q_v\): number of quadrature points per volume cell, \(e_0\): relative error at \(t=\delta t\)).}

\end{table}

\subsection{Temporal Convergence}
\label{sec:orgc399c9a}
\label{sec:temporal-convergence}
In terms of temporal convergence, we expect to have first order
in time due to use of the convex splitting scheme of \cite{gao_gradient_2012}.
We still use
the initial condition (\ref{convergence-test-init}) over a disk
domain to do perform this test. To verify the expectation, we fix
\(\epsilon = 0.5\) and the final time \(t_f = 0.1\), and perform
self-convergence tests by taking different \(\delta t\). The
numerical error is measured at \(t=t_f\) using the formula
\(e_f := \lVert \phi^f_h - \phi^f_0 \rVert_\infty / \lVert\phi^f_h\rVert_\infty\),
against the reference solution obtained by taking \(\delta t =
0.0125\).
The results are presented in Table \ref{tab:convergence-time},
confirming that the scheme is of first order in time.

\begin{table}[htbp]
\centering
\begin{tabular}{|c|c|c|}
\hline
\(\delta t\) & \(e_f\) & EOC\\
\hline
\(0.1\) & \(2.92 \times 10^{-1}\) & --\\
\(0.05\) & \(1.42 \times 10^{-1}\) & \(1.04\)\\
\(0.025\) & \(7.53 \times 10^{-2}\) & \(0.92\)\\
\hline
\end{tabular}
\caption{\label{tab:convergence-time}Temporal Convergence Tests (\(\delta t\): time step size, \(e_f\): relative error at final time).}

\end{table}

\subsection{Short-time Dynamics}
\label{sec:orgccb38cd}
\label{sec:spinodal-decomposition}
In this example, we test for the ability of our scheme to resolve small structures
using high order approximation in combination with adaptive mesh refinement.
We mimic, using the initial condition, a near-uniform state with small random perturbations.
Those perturbations will lead to a phase separation process.
We expect that our scheme is able to resolve the dynamics of structures at
all scales during the phase separation.

The setup tries to capture short-time dynamics of the Cahn-Hilliard equation.
We choose \(\epsilon = 10^{-2}\), \(\delta t = 10^{-4}\), and use fourth-order
quadrature rules for volume and boundary discretizations.
For volume densities,
\begin{equation}
\label{eqn:phi-4th-order-approximation}
  \phi(x,y) \approx \sum_{0\leq k,l\leq3} a_{k,l} P_k(x) P_l(y),
\end{equation}
(where \(P_k\) stands for Legendre polynomial of degree \(k\) defined on the
concerning cell) we use
\begin{equation}
\label{eqn:err-phi-4th-order-approximation}
  E := |a_{3,3}| + \sum_{0\leq k < 3} (|a_{k,3}| + |a_{3,k}|)
\end{equation}
as the error indicator for the cell.

The simulation starts with an initial profile generated
using the following procedure:
\begin{enumerate}
\item A uniform $200\times 200$ Cartesian grid is generated on $B$. Then
      for each grid point a real number is drawn from a uniform random
      distribution on $[-10^{-3}, 10^{-3}]$.
\item A cubic B-spline representation of the initial $\phi$ is formed
      using the Cartesian grid points as knots and the random numbers as
      interpolating values.
\item The values and derivatives of the intial $\phi$ are then evaluated
      using the cubic B-spline representation on the quadrature points
      of the box mesh as well as the boundary mesh.
\end{enumerate}

\begin{figure}
  \centering
  \subfigure[$t=0$.\label{fig:spinodal-decomposition-0}]{
    \includegraphics[width=0.3\linewidth]{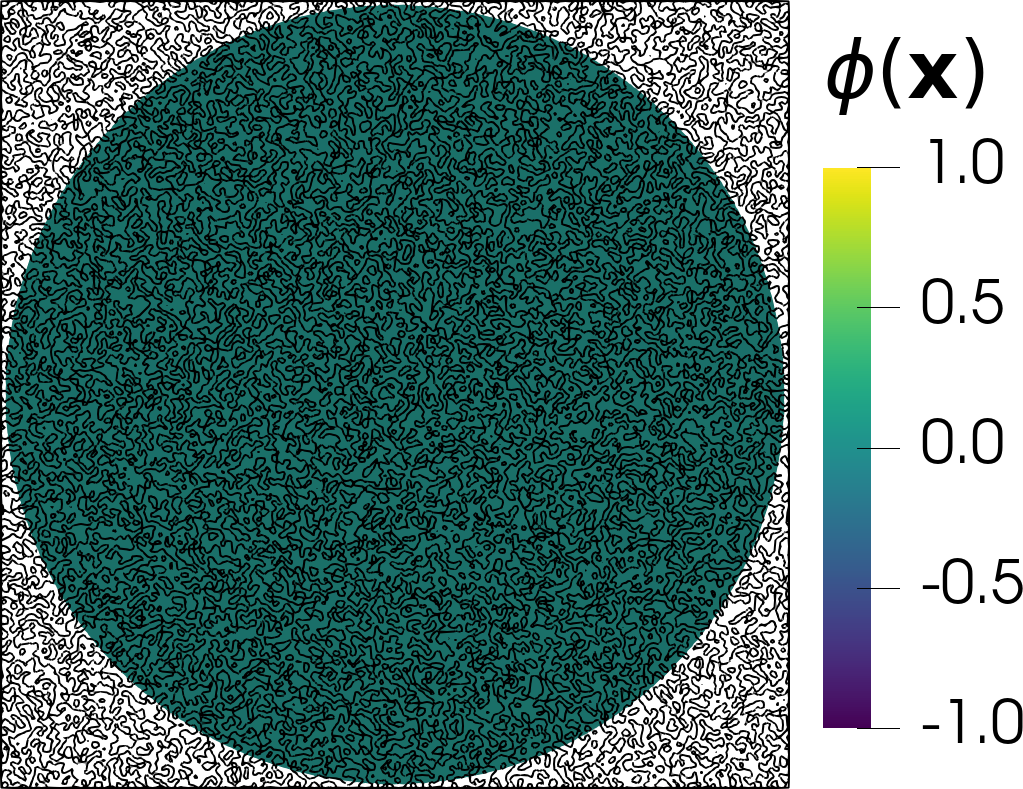}}
  \subfigure[$t=\delta t$.\label{fig:spinodal-decomposition-1}]{
    \includegraphics[width=0.3\linewidth]{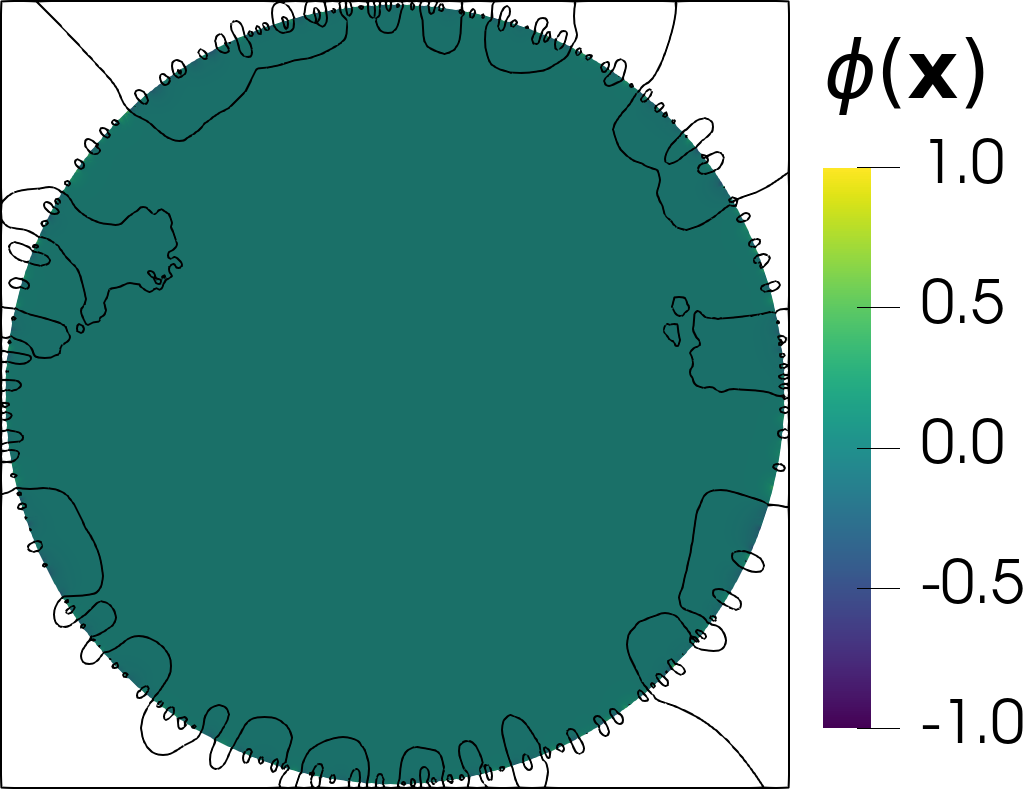}}
  \subfigure[$t=5 \delta t$.\label{fig:spinodal-decomposition-5}]{
    \includegraphics[width=0.3\linewidth]{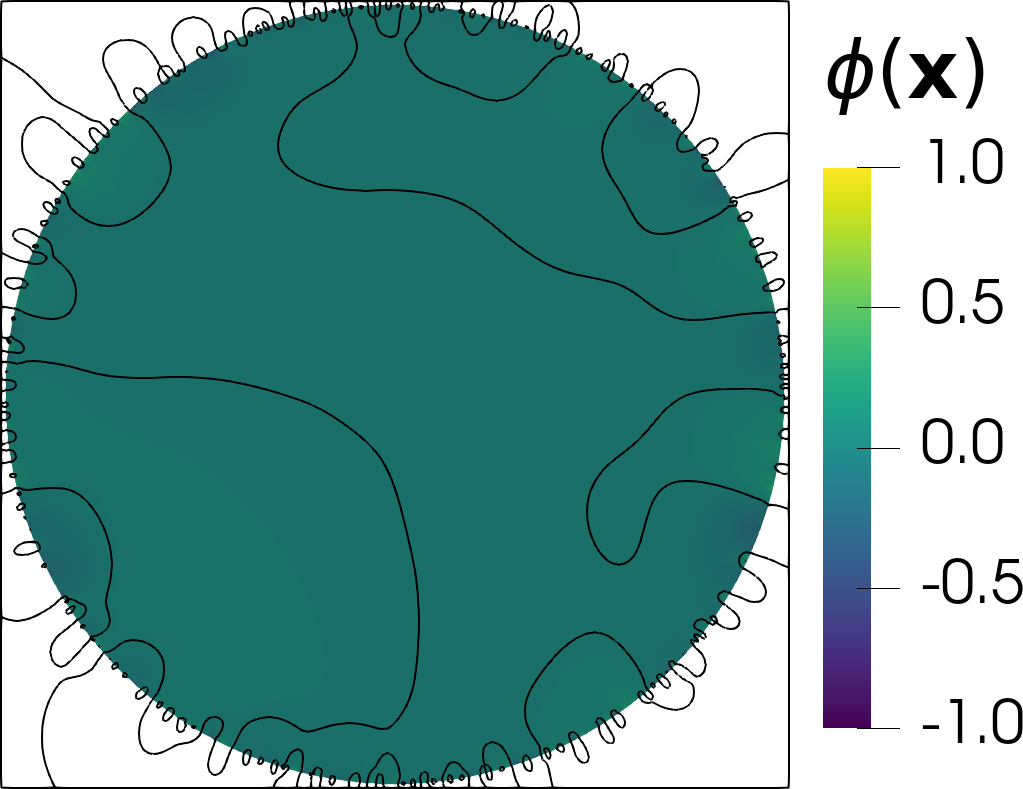}}
  \subfigure[$t=10 \delta t$.\label{fig:spinodal-decomposition-10}]{
    \includegraphics[width=0.3\linewidth]{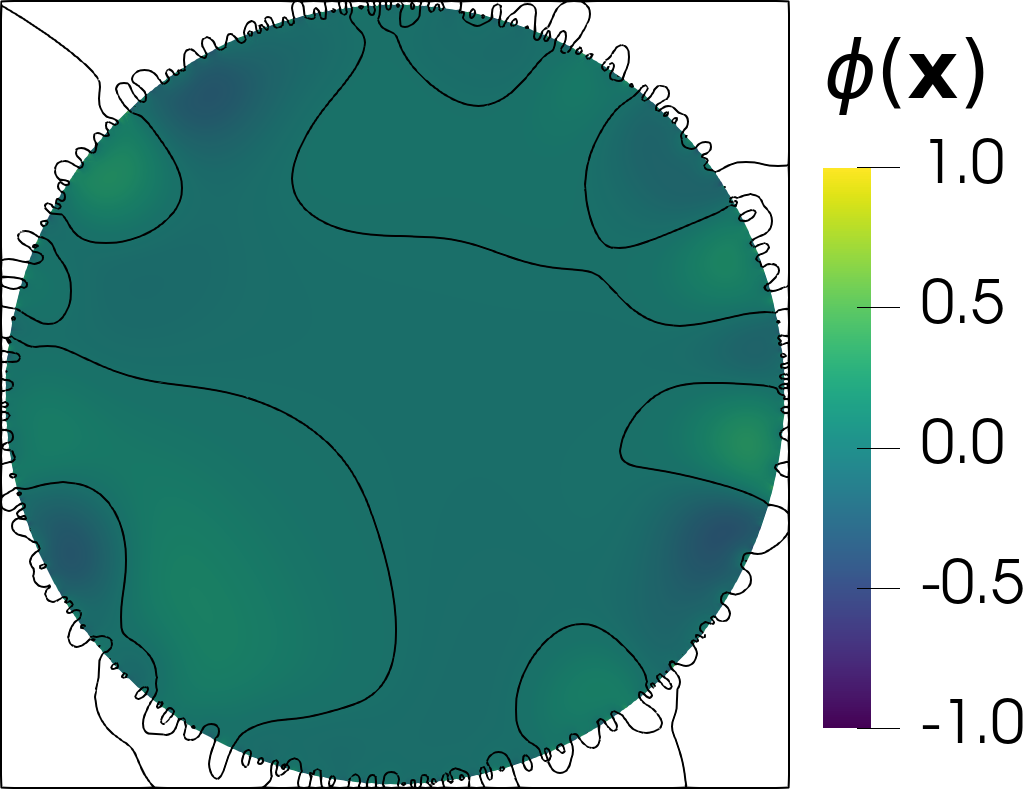}}
  \subfigure[$t=20 \delta t$.\label{fig:spinodal-decomposition-20}]{
    \includegraphics[width=0.3\linewidth]{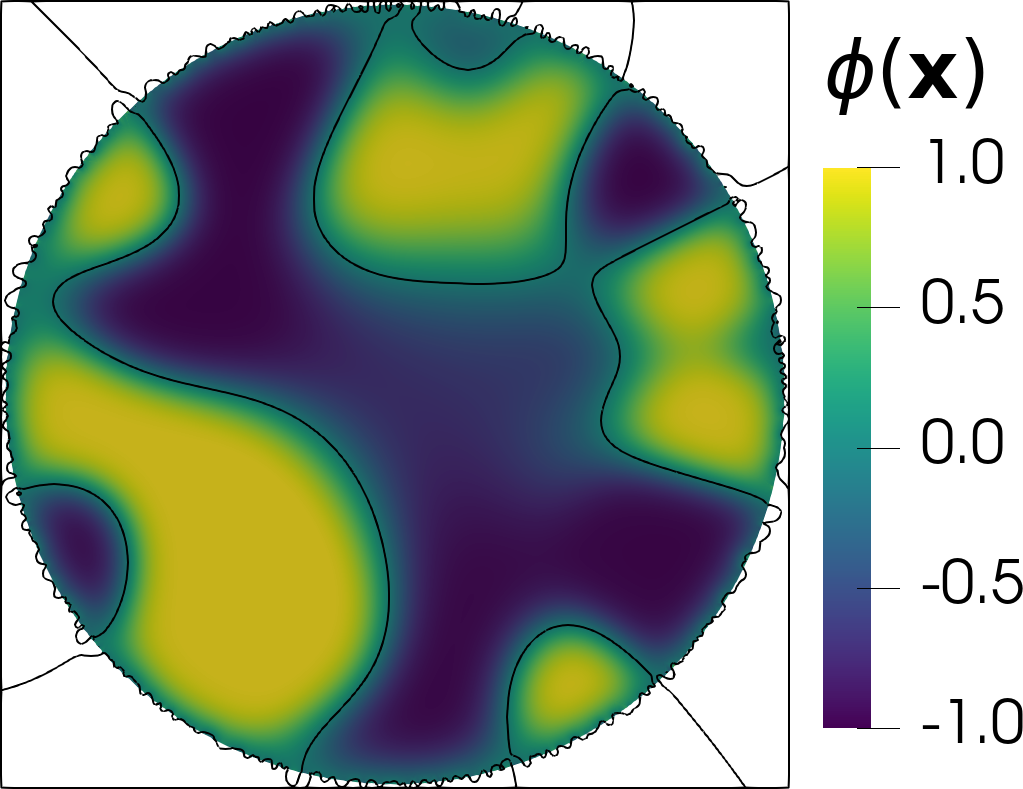}}
  \subfigure[$t=30 \delta t$.\label{fig:spinodal-decomposition-30}]{
    \includegraphics[width=0.3\linewidth]{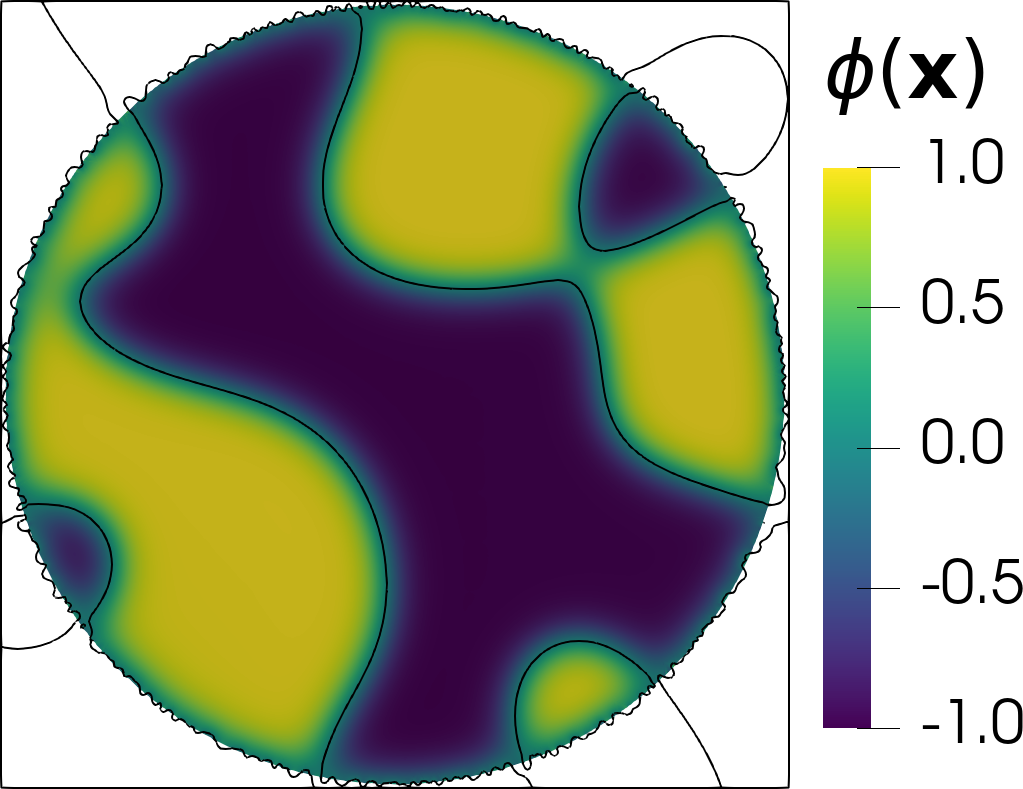}}
  \caption{Short-time Dynamics of the Cahn-Hilliard Equation ($\delta t = 10^{-4}$).}
  \label{fig:spinodal-decomposition}
\end{figure}

Numerical solutions are shown in Figure \ref{fig:spinodal-decomposition}.
The first \(3\) time steps are performed solving the full SKIE for the boundary
integral equation to capture the fast dynamics near the boundary,
while the later steps are performed using the approximation
(\ref{eqn:cahn-hilliard-representation-u-stabilized}).
In our experience, the numerical solution tends to form boundary layers
in regions away from contact lines (where its tangential derivative is close to \(0\)).
The color maps are shown for each sub-figure individually, while the contours
are all showing the zero level set of \(\phi\) over the whole bounding box.
The short-time dynamics can be understood as interplay of three mechanisms:
surface-energy driven phase transition,
bulk-energy driven spinodal decomposition, and
solid boundary induced nucleation.
In later time steps, the smoothing term in the free energy starts to play a more
important role and leads to coarsening of the patterns.
As a benefit of the high-order accuracy of our scheme, our method is able to
resolve the interfaces as well as the small patterns formed on the boundary very
well.

\subsection{Long-time Dynamics}
\label{sec:org6c056c4}
\label{sec:motion-by-curvature}
Since our solver is only first-order in time, choosing larger time step sizes
reduces solution accuracy significantly; however, through this example we
demonstrate that even with large time steps, higher spatial order is still
useful for preserving certain qualitative properties of the solution.
In this example, we verify that the solver can preserve
solution symmetry across long time scales.

We still use \(\epsilon = 10^{-2}\),
but take larger time step size \(\delta t = 0.5\).
The initial condition is
\begin{equation}
\label{eqn:long-time-initial-conditions}
\phi^0(x,y) =
  \sin \left( \frac{40\pi}{L} x \right)
  \cos \left( \frac{32\pi}{L} y \right),
\end{equation}
where \(L=0.5\) is the size of the bounding box.
Note that \(\phi^0(x,-y) = \phi^0(x,y)\).
For this example we use
(\ref{eqn:cahn-hilliard-representation-u-stabilized}) from the start.
Since \(\delta t\) is rather large compared to \(\delta x\), the temporal truncation
error of convex splitting scheme dominates; however, we expect that the
scheme is still stable and reaches an equilibrium state when
\(t\rightarrow\infty\), and that the solution \(\phi\) at any time stays
symmetric with respect to the \(x\) axis.

\begin{figure}
  \centering
  \subfigure[$t=0$.\label{fig:coarsening-0}]{
    \includegraphics[width=0.3\linewidth]{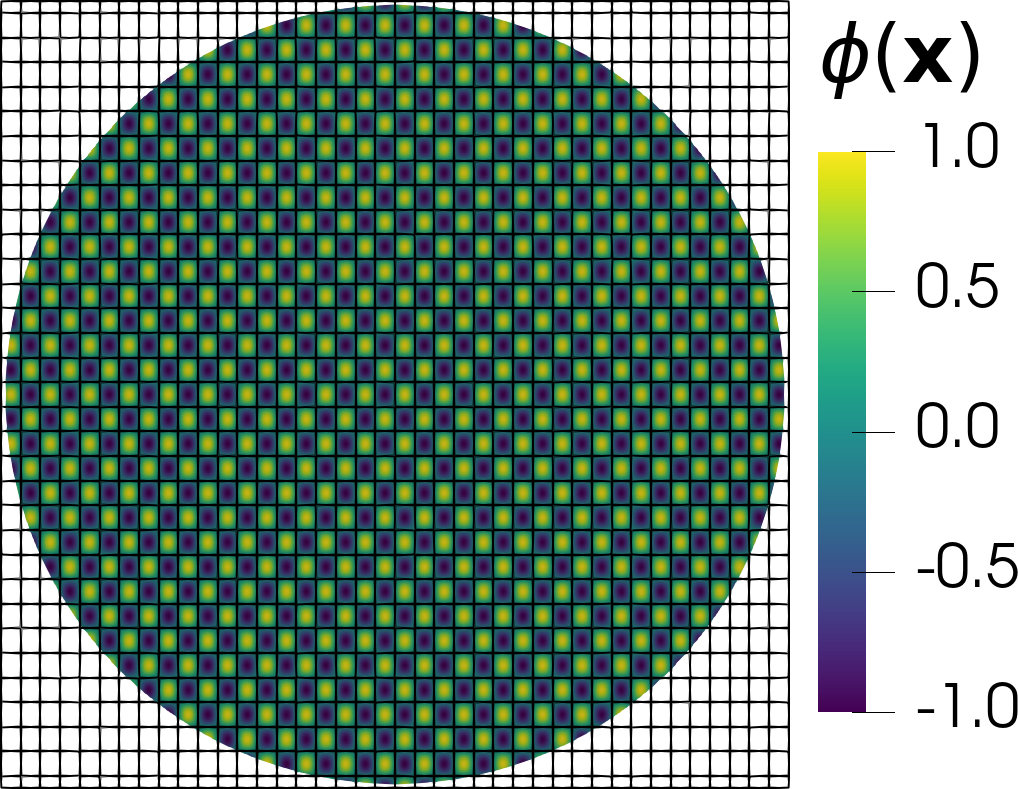}}
  \subfigure[$t=\delta t$.\label{fig:coarsening-1}]{
    \includegraphics[width=0.3\linewidth]{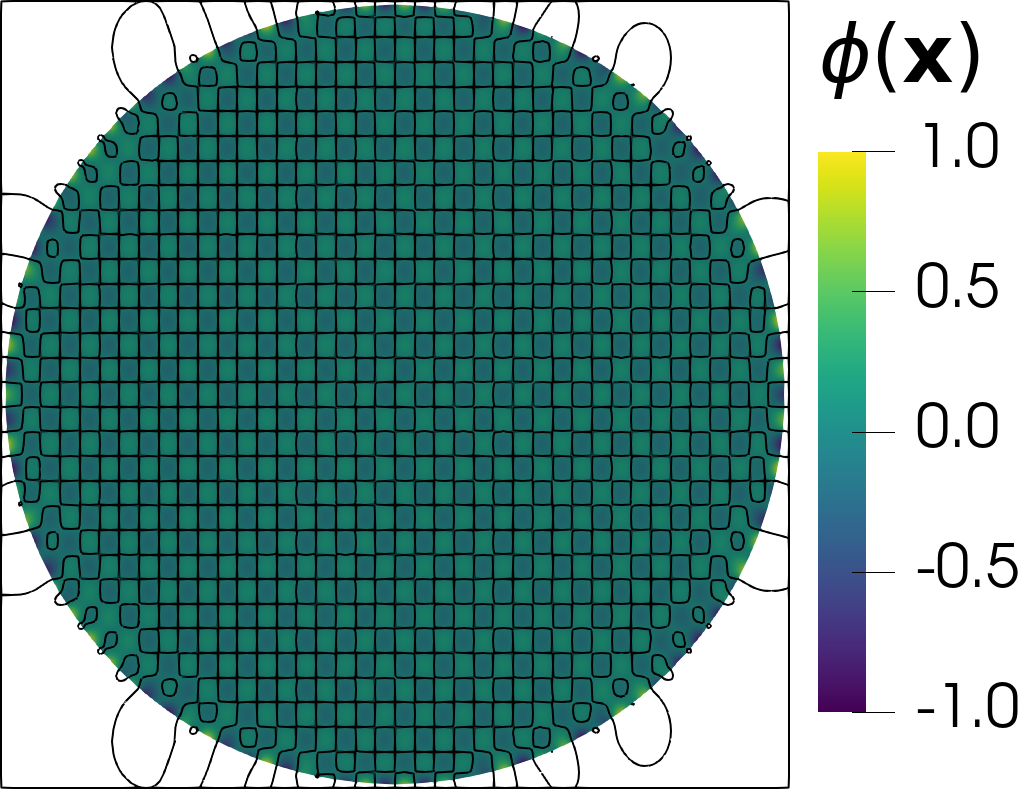}}
  \subfigure[$t=11 \delta t$.\label{fig:coarsening-11}]{
    \includegraphics[width=0.3\linewidth]{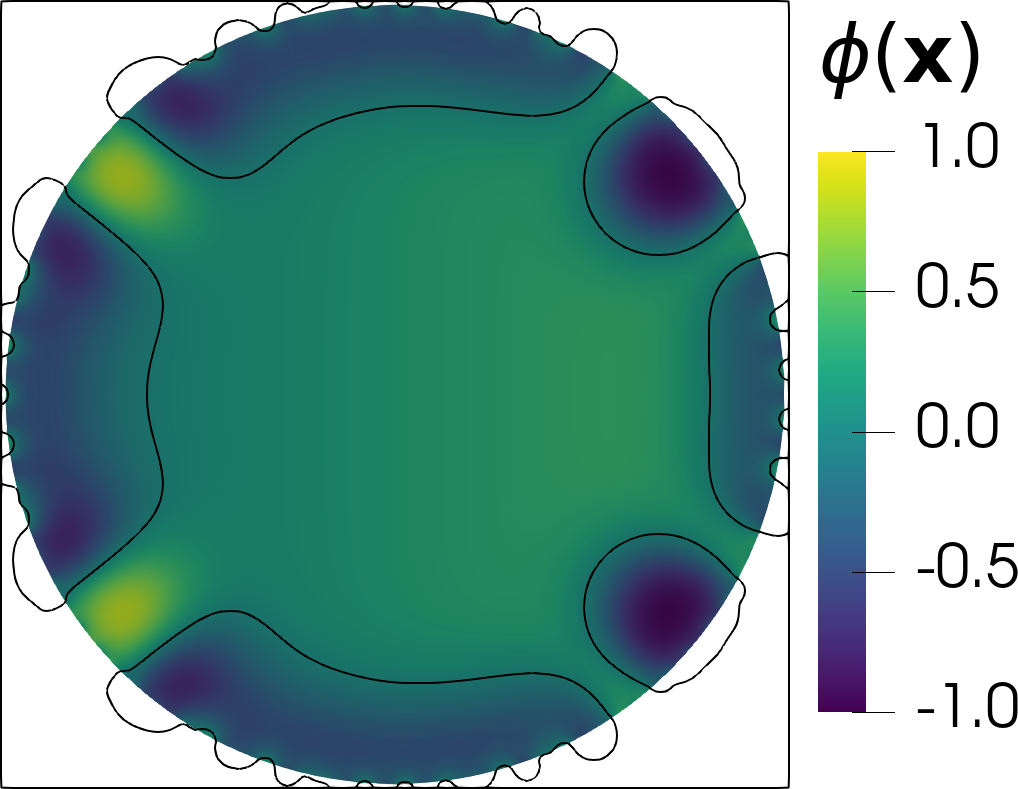}}
  \subfigure[$t=51 \delta t$.\label{fig:coarsening-51}]{
    \includegraphics[width=0.3\linewidth]{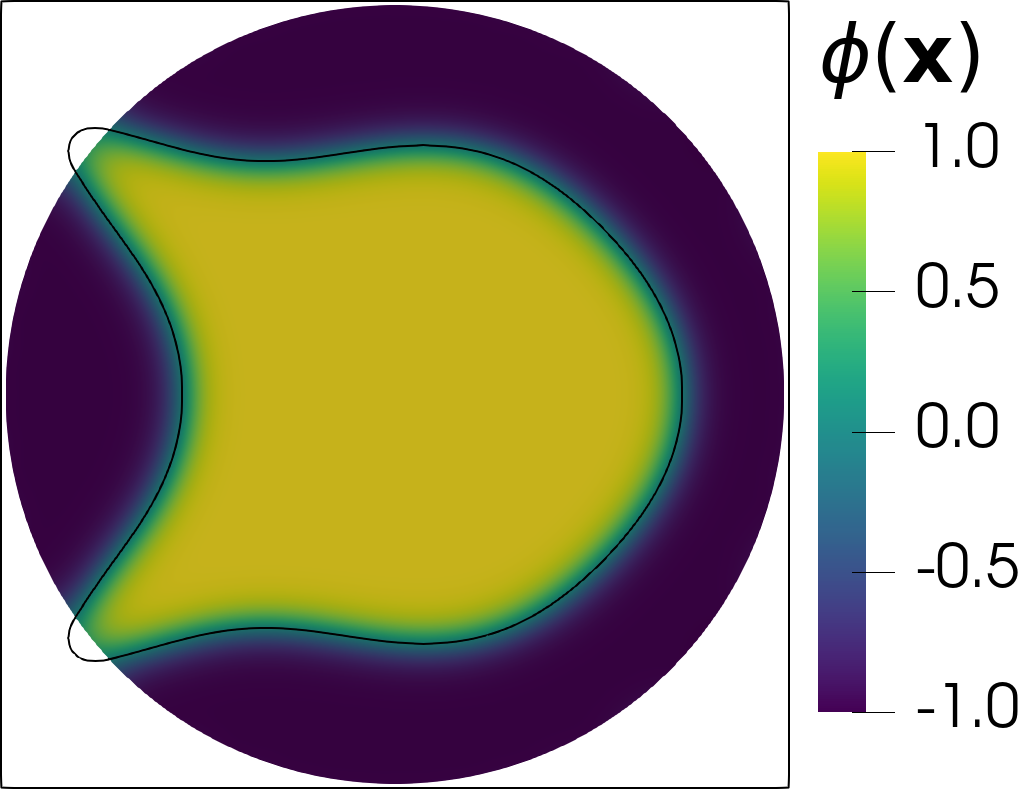}}
  \subfigure[$t=101 \delta t$.\label{fig:coarsening-101}]{
    \includegraphics[width=0.3\linewidth]{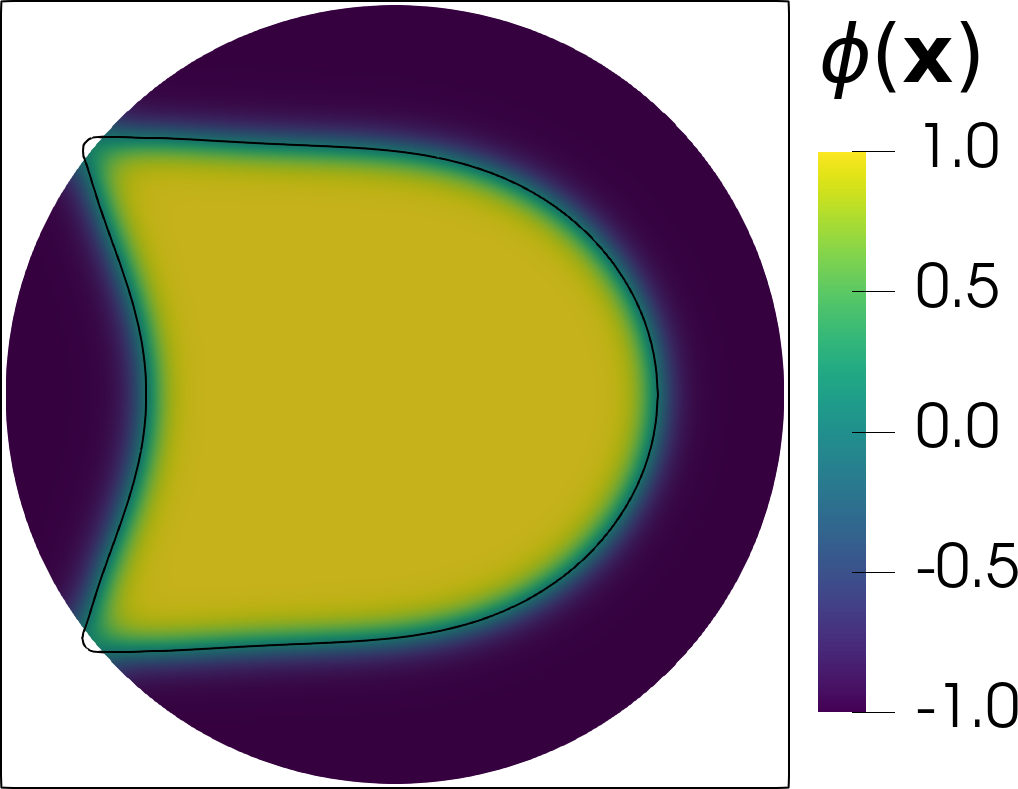}}
  \subfigure[$t=1451 \delta t$ (Equilibrium).\label{fig:coarsening-1451}]{
    \includegraphics[width=0.3\linewidth]{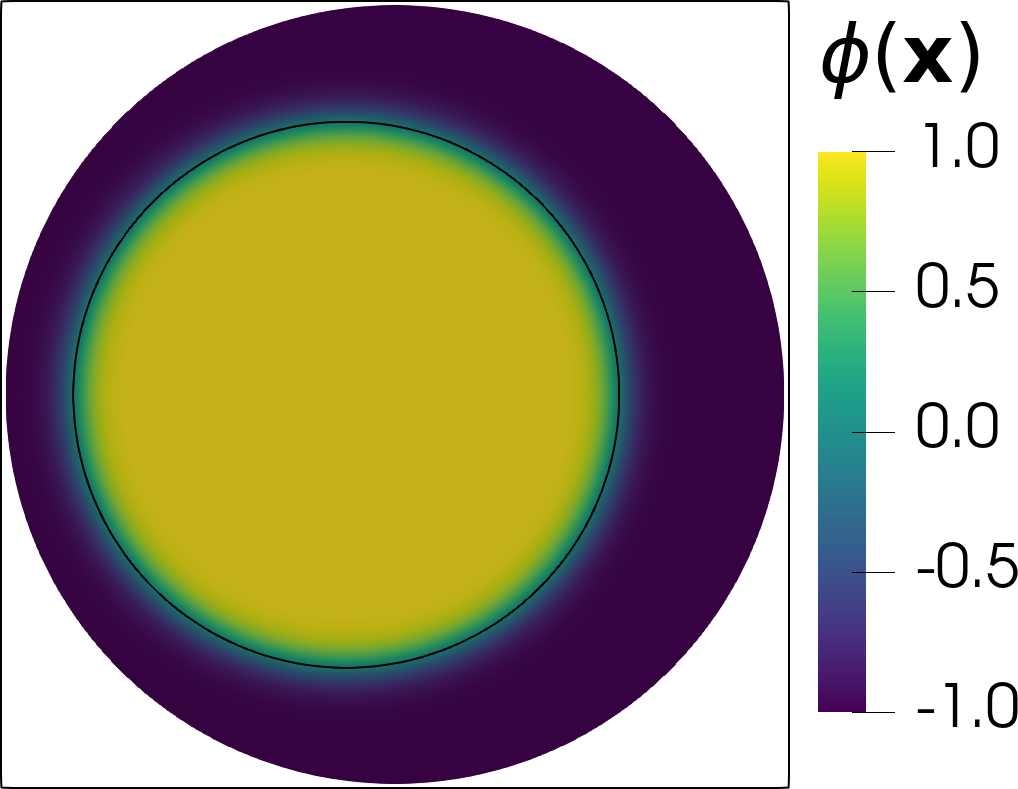}}
  \caption{Long-time Dynamics of the Cahn-Hilliard Equation ($\delta t = 0.5$).}
  \label{fig:coarsening}
\end{figure}

Numerical results are shown in Figure \ref{fig:coarsening}. The simulation stops
at \(t=780\) when \(\lVert\phi^{n+1} - \phi^n\rVert_\infty < 10^{-14}\), which
confirms that the numerical solution reaches an equilibrium state.
Besides being able to resolve the interface and the contact lines,
our method is also shown to be able to properly resolve the
interactions between the contact lines and the geometry with no
alignment requirements between the two sets of meshes:
in the case circular domain, as long as both the boundary mesh
and the box mesh are symmetric, the solution symmetry will be
well-preserved from the beginning to the equilibrium state.

\begin{figure}
  \centering
  \subfigure[Elliptical Disk ($t=0.1$).\label{fig:geometry-ell-01}]{
    \includegraphics[width=0.3\linewidth]{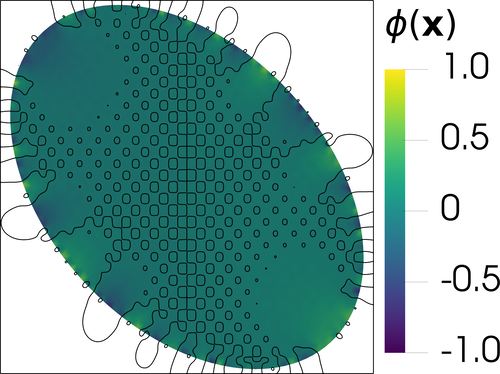}}
  \subfigure[Rounded Square ($t=0.1$).\label{fig:geometry-sq-01}]{
    \includegraphics[width=0.3\linewidth]{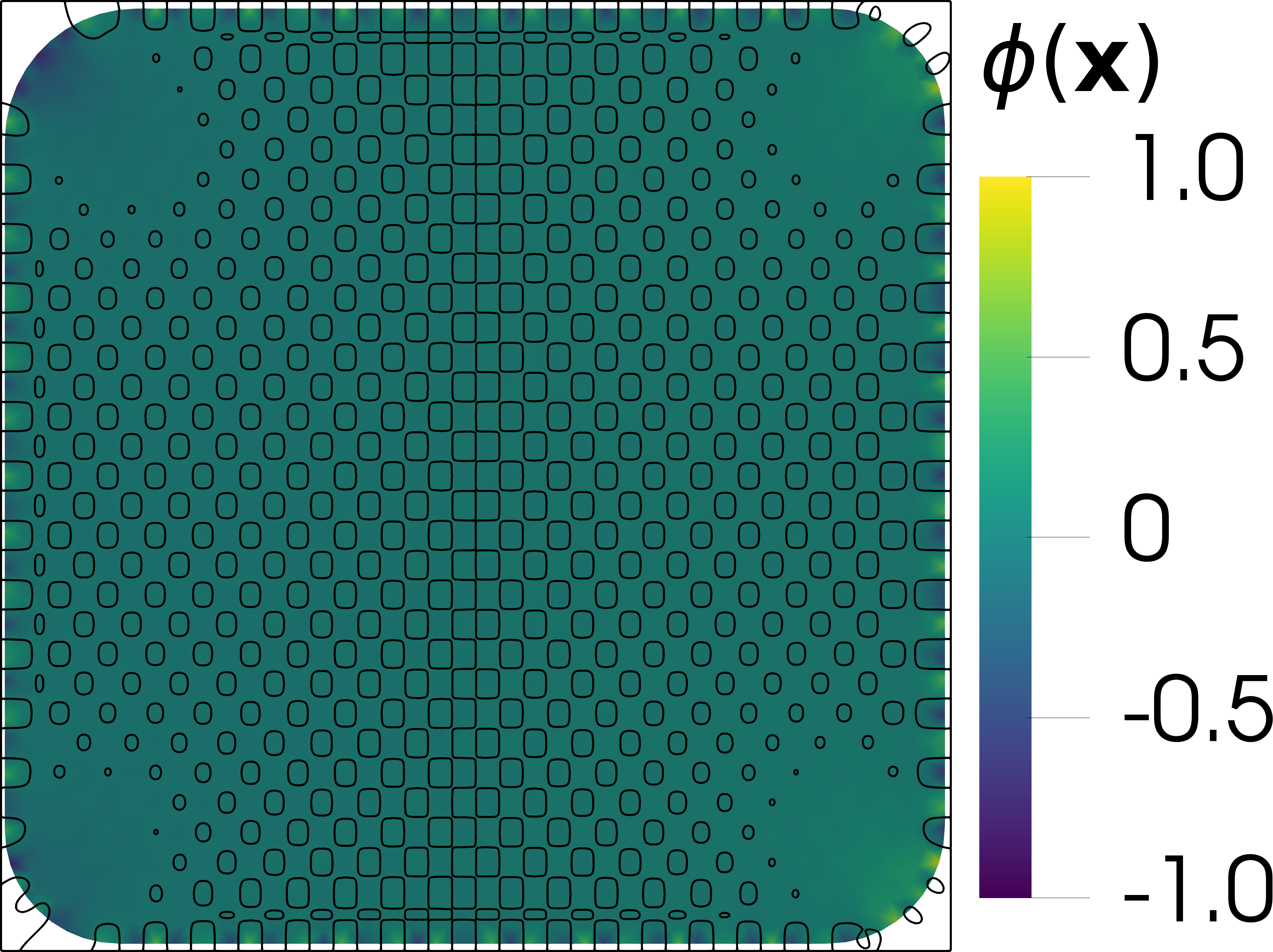}}
  \subfigure[Rounded L ($t=0.1$).\label{fig:geometry-l-01}]{
    \includegraphics[width=0.3\linewidth]{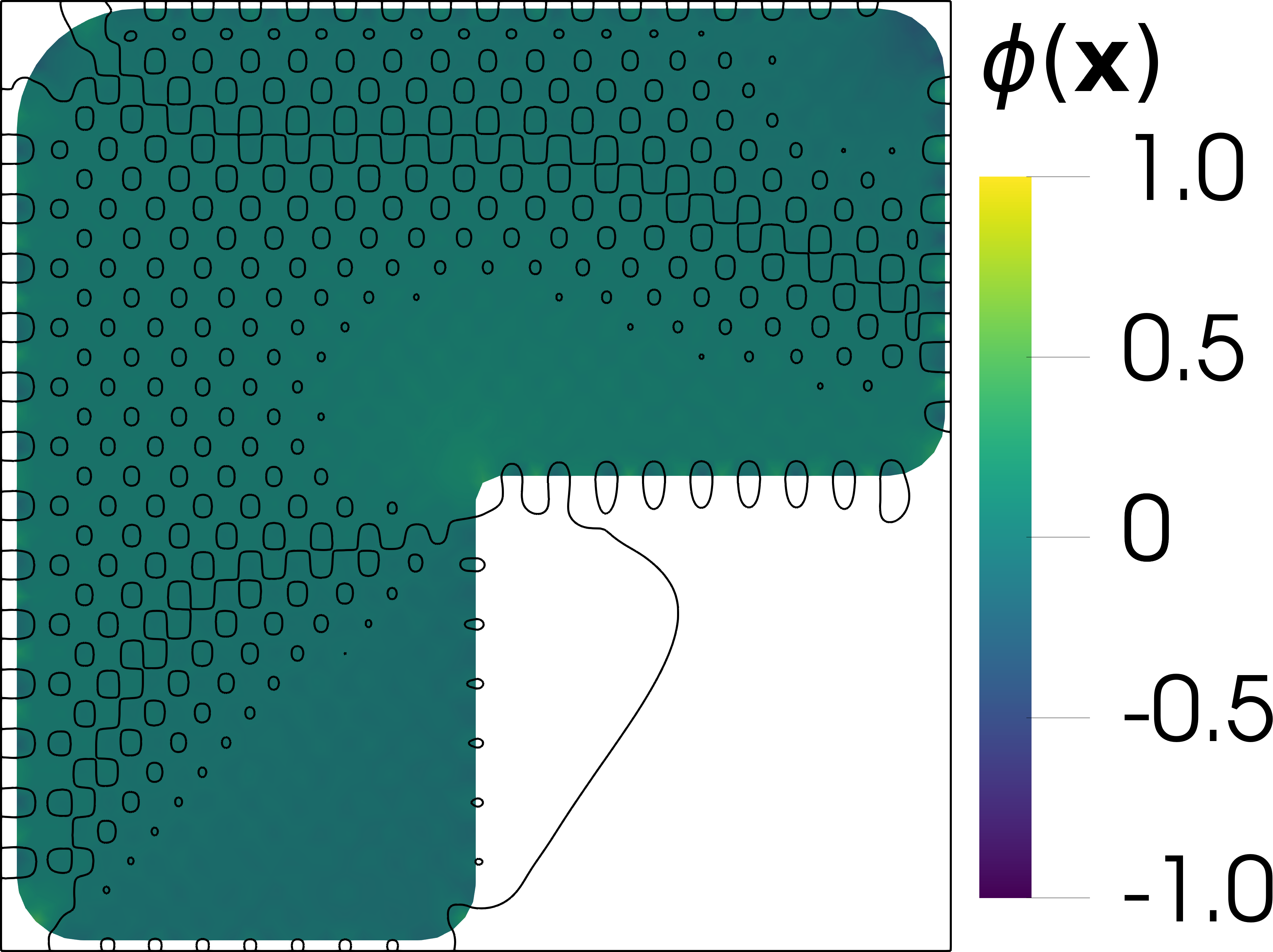}}
  \subfigure[Elliptical Disk ($t=0.5$).\label{fig:geometry-ell-05}]{
    \includegraphics[width=0.3\linewidth]{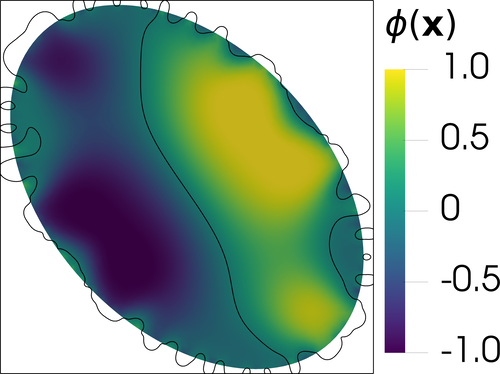}}
  \subfigure[Rounded Square ($t=0.5$).\label{fig:geometry-sq-05}]{
    \includegraphics[width=0.3\linewidth]{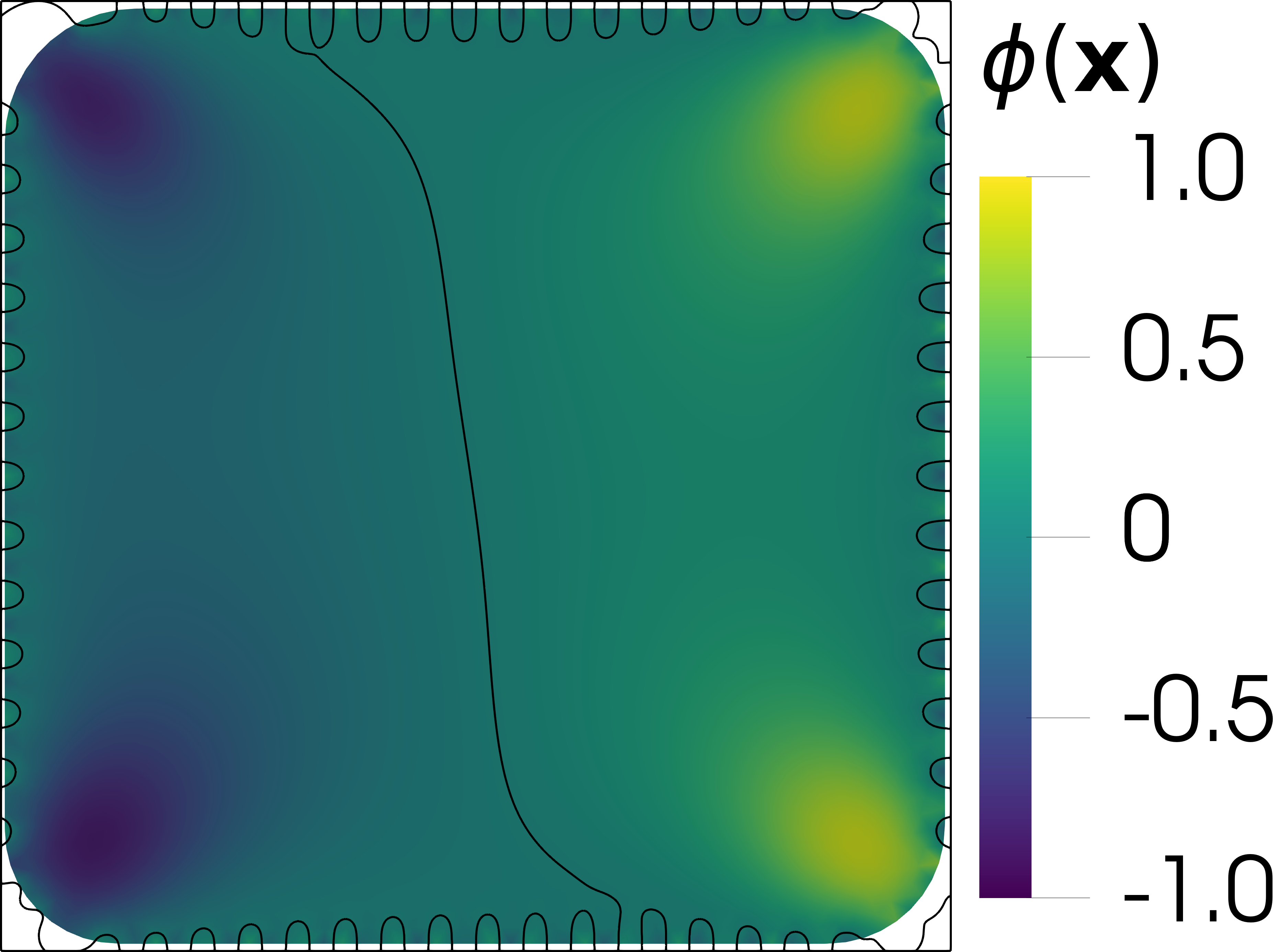}}
  \subfigure[Rounded L ($t=0.5$).\label{fig:geometry-l-05}]{
    \includegraphics[width=0.3\linewidth]{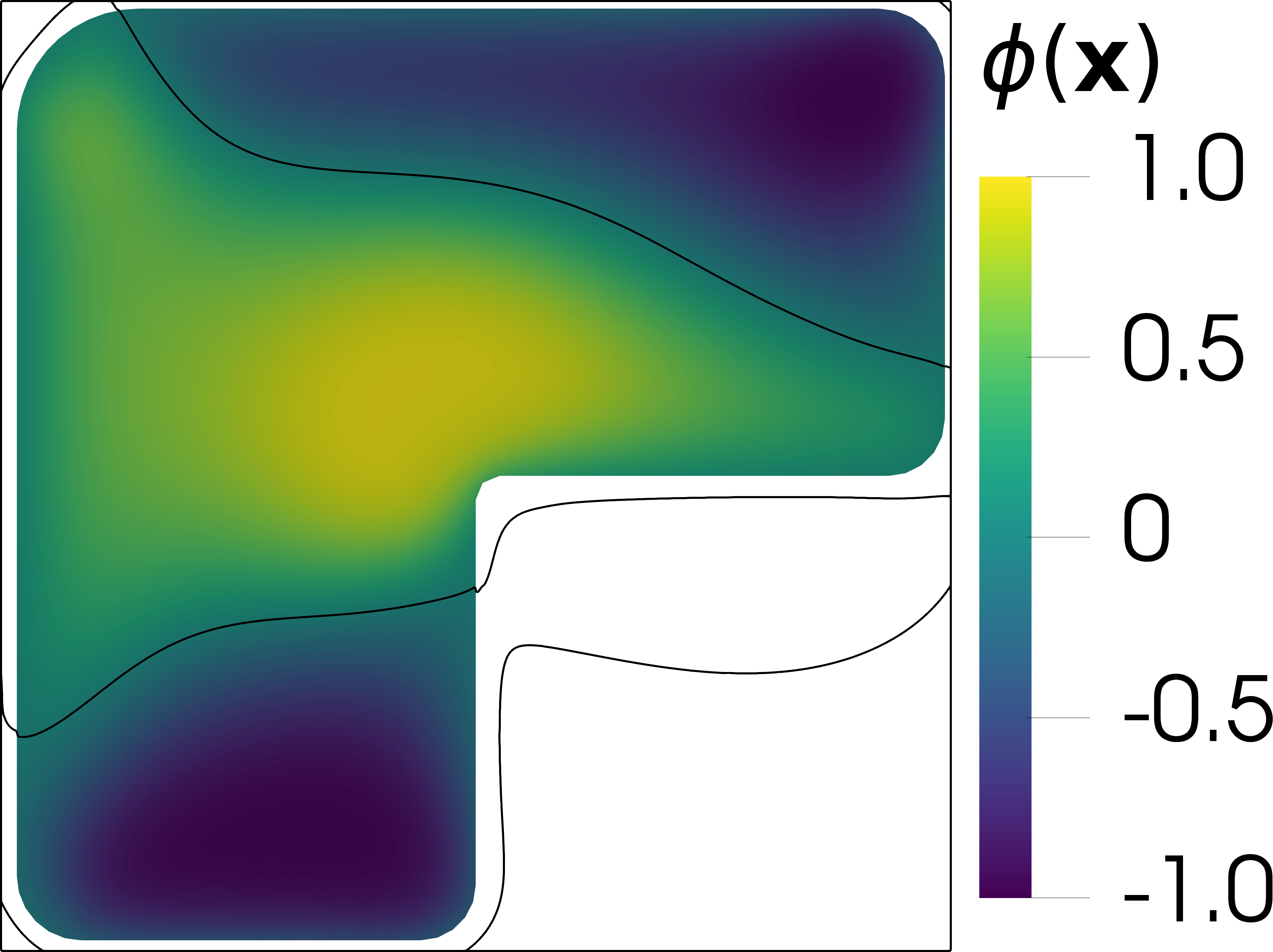}}
  \subfigure[Elliptical Disk ($t=1$).\label{fig:geometry-ell-1}]{
    \includegraphics[width=0.3\linewidth]{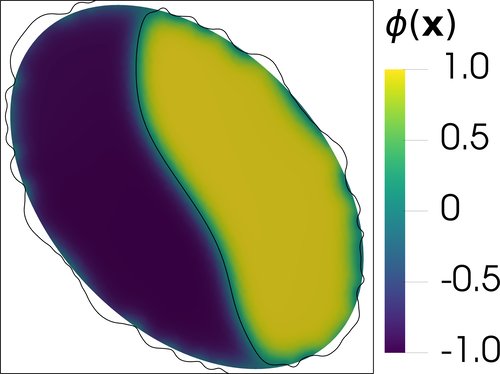}}
  \subfigure[Rounded Square ($t=1$).\label{fig:geometry-sq-1}]{
    \includegraphics[width=0.3\linewidth]{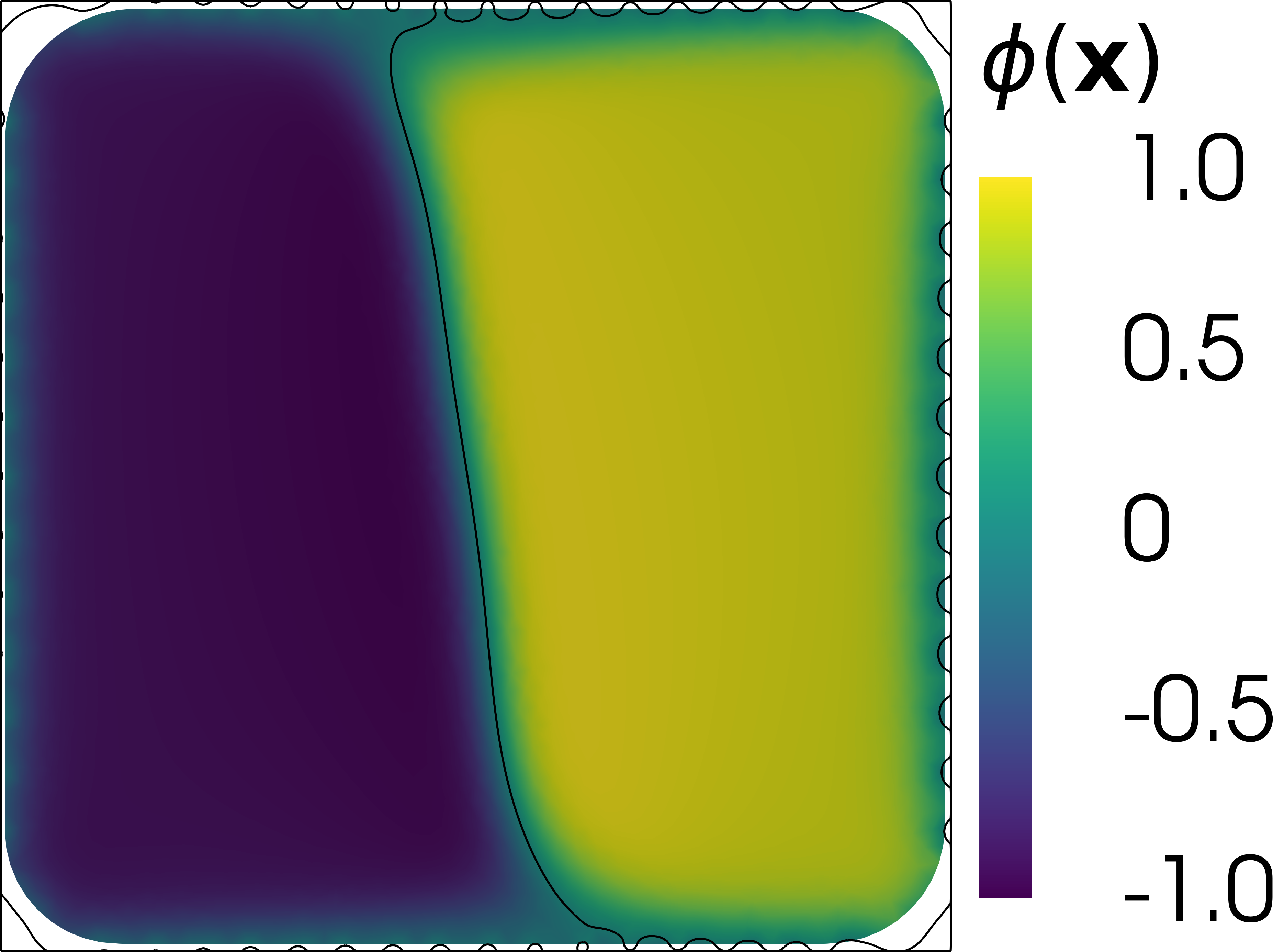}}
  \subfigure[Rounded L ($t=1$).\label{fig:geometry-l-1}]{
    \includegraphics[width=0.3\linewidth]{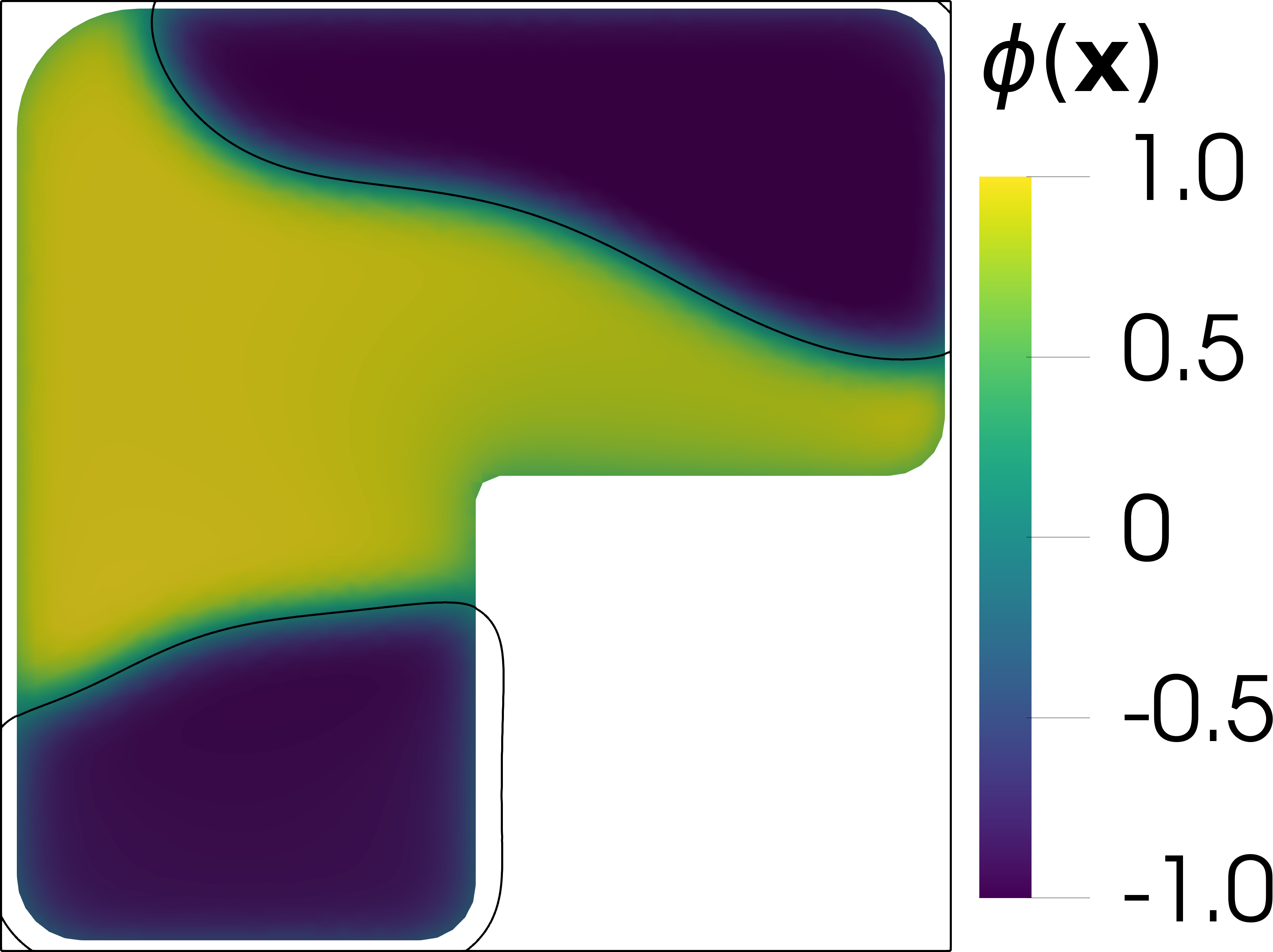}}
  \caption{Numerical Solutions for Complex Geometries.}
  \label{fig:complex-geometries}
\end{figure}

Then, we demonstrate the ability of our scheme to
handle moderately complex smooth geometries. We report the numerical
solution at \(t = 1\) for the same problem setup, but
with different physical domains. The results are shown in
Figure \ref{fig:complex-geometries}.

\begin{figure}
  \centering
  \input{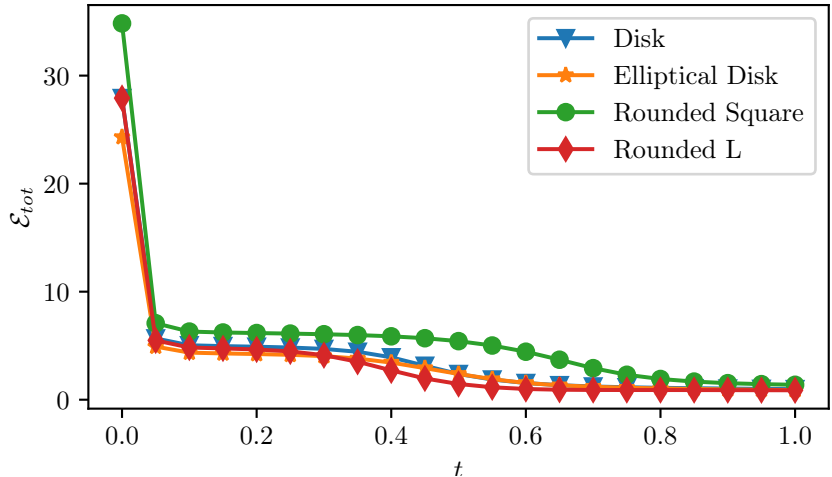}
  \caption{Total Energy of the Numerical Solutions.}
  \label{fig:energy-complex-geometries}
\end{figure}

The total energy of the numerical solution for different
geometries is shown in Figure
\ref{fig:energy-complex-geometries}. According to the theory, the
energy decaying property is preserved as long as the numerical
solution remains bounded. With the stabilization applied, we
observe that the numerical solution remains bounded for various
geometries.

\subsection{Comparisons with FEM}
\label{sec:org7a2d0a6}
\label{sec:comparisons-with-FEM}

It is a natural question to ask how the integral equation approach
that we take in this article
compares to more popular finite element methods. In this example we
aim to provide some perspective on this by presenting a comparison between
our solver and a simple finite element solver.

The finite element solver we consider is based on
a mixed formulation for the time-discretized system
(\ref{eqn:convex-splitting-phi})--(\ref{eqn:convex-splitting-bc-mu}),
using \(Q^1\) elements for both \(\phi\) and \(\mu\).
We report the condition number and number of degrees of freedom
with their dependence on \(\delta x\)
for FEM in Table \ref{tab:condition-numbers-FEM},
where \(\delta x\) is calculated as the square root of average cell measure,
and the same statistics for IEM in Table
\ref{tab:condition-numbers-IEM}. The tests are performed over a
circular domain with radius \(0.247\). In all tests with both the
FEM and the IEM, we fix \(\epsilon = 0.5\), \(\theta_Y = \pi/3\) and
\(\delta t = 1\). The rest of the parameters for the IEM are the
same as in \ref{sec:convergence-test}.
The data confirms that, being a second-kind integral equation, the
linear system has bounded condition number that does not scale with mesh
size, unlike the FEM whose linear system's condition number grows
in an unbounded manner when refining the mesh.
In addition, when \(\delta x \rightarrow 0\), the number of
unknowns in the linear solve scales as \(\mathcal{O}(1/\delta x)\) for our
method, while it scales as \(\mathcal{O}(1/\delta x^2)\) for the FEM.

\begin{table}[htbp]
\centering
\begin{tabular}{|c|c|c|}
\hline
\(\delta x\) & \(N_{DOFs}(FEM)\) & \(\kappa(FEM)\)\\
\hline
1.11\texttimes{} 10\textsuperscript{-1} & 16 & 9.79\texttimes{} 10\textsuperscript{3}\\
6.26\texttimes{} 10\textsuperscript{-2} & 50 & 1.24\texttimes{} 10\textsuperscript{4}\\
3.32\texttimes{} 10\textsuperscript{-2} & 178 & 2.15\texttimes{} 10\textsuperscript{4}\\
1.71\texttimes{} 10\textsuperscript{-2} & 674 & 6.82\texttimes{} 10\textsuperscript{4}\\
8.64\texttimes{} 10\textsuperscript{-3} & 2626 & 1.12\texttimes{} 10\textsuperscript{5}\\
\hline
\end{tabular}
\caption{\label{tab:condition-numbers-FEM}Condition Numbers of FEM Matrices (2nd Order).}

\end{table}

\begin{table}[htbp]
\centering
\begin{tabular}{|c|c|c|}
\hline
\(\delta x\) & \(N_{DOFs}(IEM)\) & \(\kappa(IEM)\)\\
\hline
1.11\texttimes{} 10\textsuperscript{-1} & 56 & 9.86\texttimes{} 10\textsuperscript{2}\\
6.26\texttimes{} 10\textsuperscript{-2} & 100 & 9.87\texttimes{} 10\textsuperscript{2}\\
3.32\texttimes{} 10\textsuperscript{-2} & 200 & 9.88\texttimes{} 10\textsuperscript{2}\\
1.71\texttimes{} 10\textsuperscript{-2} & 364 & 9.88\texttimes{} 10\textsuperscript{2}\\
8.64\texttimes{} 10\textsuperscript{-3} & 720 & 9.89\texttimes{} 10\textsuperscript{2}\\
\hline
\end{tabular}
\caption{\label{tab:condition-numbers-IEM}Condition Numbers of IEM Matrices (4th Order).}

\end{table}

\section{Conclusions}
\label{sec:org4c552a4}
\label{sec:conclusions}

We have introduced an integral equation method for the
Cahn-Hilliard equation in bounded 2D domains with
solid boundary conditions.
For each time step, the method consists of two stages:
(1) evaluate volume potentials and
(2) solve remaining SKIE.
Fast multipole based fast algorithms are used to attain linear complexity for
both stages. To handle complex geometry, a \(C^1\) extension
of \(\phi\) onto a larger domain is performed by solving an additional SKIE.
Through numerical experiments, the method is shown to require fewer
degrees of freedom and maintain bounded condition number compared to finite element
method. To remedy some of the numerical stability concerns for
long-time simulations caused by boundary layers, we have proposed a stabilized
representation approximating the original problem, and we have validated the
approximation with numerical results.

The present work opens some exciting future perspectives.
We will seek to apply the same techniques in this paper to second-order time
discretization schemes, for example the SAV method
\cite{shen_scalar_2018} offering
higher accuracy order in time as well as superior energy stability.
In addition, the challenges posed by boundary layers are by no means unique to
Cahn-Hilliard equations, and it is of great interest to develop fast volume
potential evaluators that remain stable and accurate when the source density
develops boundary layers.
Once the boundary layer issue can be handled, our SKIE
formulation will be immediately applicable to the full long-time simulation.

\section{Acknowledgments}
\label{sec:orgcad8af3}

This research was supported in part by the National Science Foundation under
grants DMS-1418961 and DMS-1654756, and by the Hong Kong Research Grant Council
(RGC-GRF grants 605513 and 16302715, RGC-CRF grant C6004-14G, and NSFC-RGC joint
research grant N-HKUST620/15).
Any   opinions,   findings,   and conclusions, or   recommendations expressed in
this article are those of the authors and do not necessarily reflect the views
of the National Science Foundation or the Hong Kong Research Grant Council;
Neither NSF nor HKRGC has approved or endorsed its content.
Part of the work was performed while the authors were participating in
the HKUST-ICERM workshop `Integral Equation Methods, Fast
Algorithms and Their Applications to Fluid Dynamics and Materials
Science’ held in 2017.
\label{sec:org28caa3e}
\bibliography{../references/refs}
\end{document}